\documentclass{article}
\usepackage[top=1.0in, left=1.2in, right=1.2in, bottom=1.0in]{geometry}
\usepackage[numbers,sort&compress]{natbib}
\usepackage{graphicx}
\usepackage{amsmath}
\usepackage{dsfont}
\usepackage{amssymb}
\usepackage{setspace}
\usepackage{textgreek}
\usepackage{titlesec}
\usepackage{caption}
\usepackage{xcolor}
\usepackage[most]{tcolorbox}
\usepackage[colorlinks=true,
            linkcolor=blue,
            citecolor=blue,
            urlcolor=blue,
            pdfborder={0 0 0}]{hyperref}
\usepackage[version=4,arrows=pgf-filled,textfontname=sffamily,mathfontname=mathsf]{mhchem}

\usepackage{setspace}
\usepackage{amsthm}
\theoremstyle{plain}

\newtheorem{proposition}{Proposition}[section]

\newtheorem{lemma}{Lemma}[section]
\theoremstyle{definition}
\newtheorem{remark}{Remark}[section]

\newtheorem{definition}{Definition}[section]
\usepackage{authblk}
\renewcommand{\vec}[1]{\boldsymbol{#1}}

\title{Structural functional identifiability of unknown \\ differential equation observation maps}
\author[1]{Torkel E Loman}
\author[2]{Alexander P Browning}
\author[1]{Ruth E Baker}
\date{}

\affil[1]{Mathematical Institute, University of Oxford, Oxford, United Kingdom}
\affil[2]{School of Mathematics and Statistics, University of Melbourne, Melbourne, Australia}

\begin{document}

\maketitle

\begin{abstract}

In many differential equation models of real-world systems, the mechanistic structure of the dynamics is known, but the map from the dynamics to the observations is not. This is potentially problematic, as it can prevent us from e.g. learning model parameter values from data. Examples where this occurs include measuring a fluorescent marker as a proxy for the dynamics of a cellular circuit, and using wastewater measurements to follow the development of an epidemic. The nascent field of \emph{hybrid modelling}, in which machine learning components (like neural networks) are embedded within mechanistic models, provides the tools to learn such maps. Here, the system's dynamics keep their mechanistic form, while a neural network represents the unknown observation map and is trained alongside any unknown core system parameters. This raises a question of \emph{identifiability}: when can the true observation map be learnt, and when does one simply learn some map that transforms the dynamics into the available data? Building on previous work, we show how \emph{structural functional identifiability analysis} can establish whether a model's unknown observation maps and core system parameters are \emph{identifiable}, i.e. can be uniquely learnt from data. We both prove general results on when models with unknown observation maps are, and are not, identifiable, and demonstrate a general methodology (based on \emph{differential algebra}) for assessing the identifiability of such models. This paves the way for applying classic differential equation workflows to the many real-world systems that contain such unknown maps.

\end{abstract}


\section{Introduction} \label{section:introduction}

Differential equation models are widely applied to e.g. generate predictions of real-world systems. Across many applications, however, to do so one must first learn the values of any unknown model parameters. This is typically done by collecting system data, to which the model can be fitted. Here, one can rarely observe system dynamics completely and directly. Instead, a subset of system states are observed, and often only using proxies for the real state. One example is a biochemical circuit observed through a fluorescent marker, but where the exact relationship between the activity of a circuit component and the fluorescent marker is unknown~\cite{de_jong_experimental_2010,ronen_assigning_2002}. Another is an animal population monitored indirectly through markers such as faecal remains rather than through direct observations of individuals~\cite{marques_estimating_2001}. In both cases, the system can be modelled by a differential equation whose form is known (up to the values of its parameters), but where the function mapping the system state to the observations is unknown. 

This generates \emph{differential equations with unknown observation maps}. These models contain two types of unknowns: (i) unknown scalar parameters, and (ii) one (or more) unknown functions describing how the core system dynamics map to the observations. Differential equations that combine unknown parameters and unknown functions in the core dynamics have been introduced under the term \emph{universal differential equations} (UDEs)~\cite{rackauckas_universal_2021}. UDEs are used to model systems where the governing dynamics are not fully known. However, they typically assume that the mapping between the core system state and the observations is known. The case considered in this work, where this mapping is unknown, can be considered a special case of the UDEs~\cite{philipps_current_2025}.

\begin{figure}[!htb]
    \centering
    \includegraphics[width=0.98\linewidth]{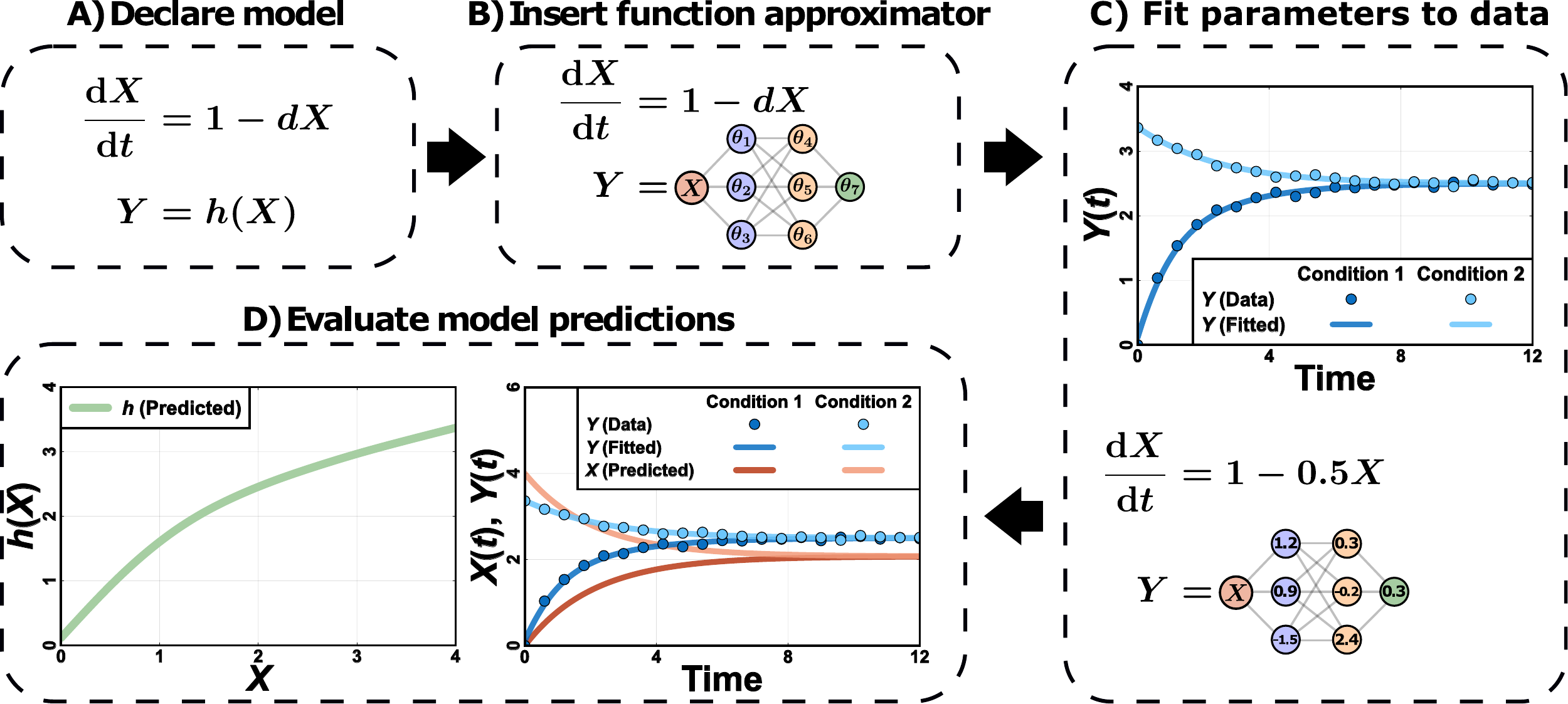}
    \caption{\textbf{Workflow for fitting a differential equation with an unknown observation map.} (A) A model is declared, here $\dot{X}=1-dX$ with the observed output $Y=h(X)$. The decay rate $d$ and the observation map $h$ are unknown. (B) The unknown map $h$ is replaced by a function approximator, here a neural network parameterised by the scalar parameters $\theta_1,\ldots,\theta_7$. (C) The upper part shows two data sets, each consisting of measurements of the observed variable $Y$. The model in B is fitted so that its $Y$ trajectories match the data points. The fitting determines the values of the unknown mechanistic and neural network parameters (lower part). (D) The right panel shows the functional form of the observation map $h$ that has been recovered from the data (as the fitted neural network). The left panel shows the data sets and fits from C, together with the recovered trajectories of the unobserved variable $X$ (which are predicted by applying the fitted $h^{-1}$ to the fitted $Y$ trajectories).}
    \label{fig:learning_obs_map_workflow}
\end{figure}

In practice, UDEs use universal function approximators, such as neural networks, to learn unknown functions from data (Figure~\ref{fig:learning_obs_map_workflow}). This creates a so-called \emph{hybrid modelling} approach, where the model combines mechanistic modelling (the model equations) with data-driven approaches (the neural network used to learn the unknown dynamics, i.e. functions). A wide range of hybrid modelling approaches have been developed~\cite{brunton_discovering_2016,chen_neural_2018,lagergren_biologically-informed_2020} and see increasingly widespread applications~\cite{bonnaffe_neural_2021,lu_neural-ode_2021,quaghebeur_hybrid_2022,kraft_towards_2022,rodriguez_einns_2023,noordijk_rise_2024,de_rooij_conditional_2025,korsbo_deepnlme_2026}. While a model could have unknown functions both in the core dynamics, and in the observation map, in this work we only consider the case where the unknown functions are in the observation map. Indeed, specialised numerical methods have been developed for this model category, and these methods are computationally more efficient than those applicable to models with unknown functions in the core dynamics~\cite{persson_petab_sciml_2026}.

Although models with unknown observation maps have wide potential applicability, they raise the question of \emph{identifiability}: to what extent is it possible to learn the true map from data, and to what extent does one simply learn an arbitrary map that transforms the dynamics into the available data? Next, under what conditions can we then learn values of core system parameters, and how is this affected by the identifiability of the observation map?

Identifiability analysis is an established concept in the literature on fitting scalar parameters. A parameter is said to be \emph{identifiable} if it can be uniquely determined from data~\cite{bellman_structural_1970,raue_structural_2009,simpson_parameter_2026}. Identifiability is generally divided into \emph{structural} and \emph{practical} identifiability. Structural identifiability is a property of the proposed model and the measurable outputs, and asks whether it is fundamentally possible to uniquely identify unknowns given an infinite quantity of noise-free data~\cite{bellman_structural_1970,cobelli_parameter_1980,wieland_structural_2021}. Practical identifiability, on the other hand, is a property of the data itself, and asks whether it is possible to recover unknowns given a finite quantity of noisy data~\cite{raue_structural_2009,wieland_structural_2021,simpson_profile-wise_2023}. 

Here, we focus on structural identifiability. In our previous work, we demonstrated how the principles of identifiability for scalar parameters could be generalised to \emph{functional identifiability} for unknown functions embedded within the right-hand sides of differential equations~\cite{loman_struct_func_2026}. In this paper, we extend this approach to where the observation maps are unknown functions. Specifically, we demonstrate how structural identifiability analysis can be applied to differential equations with unknown observation maps. We provide a general framework, based on \emph{differential algebra}~\cite{ljung_global_1994,bellu_daisy_2007}, for determining whether such models are identifiable---that is, whether the observation map can be uniquely recovered from a specific type of data. We also determine the identifiability of the core system parameters, which, notably, are sometimes identifiable even when the observation map is not. We demonstrate this framework using practically relevant examples, including models of epidemics and predator--prey systems. Finally, we provide a range of more general results for systems with unknown observation maps. This includes systems, like those involving fluorescent markers, where the unknown map has additional known structure.

\section{Mathematical framework} \label{section:methods}

\subsection{Differential equations with unknown observation maps} \label{section:obs_map_equation_def}

We will primarily consider ordinary differential equation models of the form
\begin{align}
    \dot{\vec{X}} &= \vec{F}(\vec{X};\, \vec{p}), \label{eq:base_model_ode} \\
    \vec{Y}       &= \vec{h}(\vec{X}),            \label{eq:base_model_obs}
\end{align}
where $\vec{X}\in\mathbb{R}^n$ is the state vector, $\vec{F}$ is a known functional form depending on the unknown parameters $\vec{p}\in\mathbb{R}^l$. This is the core system, which is observed through $\vec{h}=(h_1,\ldots,h_k)$, a collection of unknown, smooth \emph{observation maps} giving the observed outputs $\vec{Y}\in\mathbb{R}^k$. Each $h_j$ has a scalar output and depends on a known combination of the state variables $X_i$. While alternative cases could be considered, throughout this work we assume that the initial conditions $\vec{Y}(0)$ are known while $\vec{X}(0)$ are unknown.

In Section~\ref{section:reporters} we will consider models which are not of the form in Equations~\eqref{eq:base_model_ode}--\eqref{eq:base_model_obs}. However, for the purpose of identifiability analysis, they can be mapped to models of this form. This enables us to reuse several of our results for unknown function map-type models (as defined in Equations~\eqref{eq:base_model_ode}--\eqref{eq:base_model_obs}) to more generally defined model classes.

\subsection{Structural functional identifiability definition}
In this work, we consider the structural functional identifiability of the unknowns $\vec{h}$ and $\vec{p}$ for models of the form in Equations~\eqref{eq:base_model_ode}--\eqref{eq:base_model_obs}. Specifically, we use the following definition of structural identifiability:

\begin{definition} \label{def:functional_identifiability}
For a model of the form in Equations~\eqref{eq:base_model_ode}--\eqref{eq:base_model_obs}, consider two choices of the unknowns, $(\vec{h}_1,\vec{p}_1)$ and $(\vec{h}_2,\vec{p}_2)$. Assume that these give the same observed outputs $\vec{Y}$ over the full set of admissible initial conditions\footnote{Write $\vec{Y}(t;\vec{X}_0,\vec{h},\vec{p})$ for the output, at time $t\geq0$, of the model with unknowns $(\vec{h},\vec{p})$ and initial condition $\vec{X}(0)=\vec{X}_0$, and let $\Omega$ be the set of initial conditions consistent with prior knowledge. The two choices give the same observed outputs if $\{\vec{Y}(\cdot\,;\vec{X}_0,\vec{h}_1,\vec{p}_1):\vec{X}_0\in\Omega\}=\{\vec{Y}(\cdot\,;\vec{X}_0,\vec{h}_2,\vec{p}_2):\vec{X}_0\in\Omega\}$, i.e. if every output trajectory of one choice is also an output trajectory of the other. The two trajectories may come from different initial conditions.}. Then, an unknown function $h_j$ is:
\begin{itemize}
    \item \emph{Globally structurally identifiable} if this assumption always implies $h_{j,1}=h_{j,2}$.
    \item \emph{Locally structurally identifiable} if, for any $(\vec{h}_1,\vec{p}_1)$, there exists $\varepsilon>0$ such that this assumption implies $h_{j,1}=h_{j,2}$ whenever $\sup_{\vec{x}}|h_{j,1}(\vec{x})-h_{j,2}(\vec{x})|<\varepsilon$ holds\footnote{This is equivalent to there being a \emph{countable} number of functions $h_j$ giving the same observed dynamics. For systems with rational right-hand sides, a \emph{finite} number of functions suffice.}.
    \item \emph{Structurally non-identifiable} otherwise.
\end{itemize}
\end{definition}
 Generally, when the term ``structurally identifiable'' is used without qualification, it refers to global structural identifiability. A similar definition applies to the unknown parameters in $\vec{p}$.


\subsection{Differential algebra} \label{section:differential_algebra}

In Sections~\ref{section:diff_algebra} and \ref{section:reporters}, we will use differential algebra to determine whether models are identifiable. This approach first removes the unobserved state variables from an ODE model through repeated differentiation and substitution~\cite{ljung_global_1994,loman_struct_func_2026}. What remains is one (or multiple) \emph{input--output equation}, a differential equation relating the observed variables, and their derivatives, to the model's unknown parameters and functions. Two choices of the unknowns that give the same input--output equation cannot be distinguished from the observed output, and the existence of such choices therefore implies structural non-identifiability.

Structural identifiability considers the case of arbitrary amounts of noiseless data, which means that the observed variables and their derivatives, up to one less than the order of the input--output equation, can be chosen independently through the initial conditions. In practice, this allows each term of the input--output equation to be considered in isolation, from which the identifiable quantities can be recovered.

As an example, consider structural parameter identifiability for the model
\begin{align*}
    \dot{X}_1 &= -aX_1 + bX_2, \\
    \dot{X}_2 &= cX_1 - dX_2,
\end{align*}
with only $X_1$ observed and where $a$, $b$, $c$, and $d$ are unknown parameters. Eliminating the unobserved variable $X_2$ by differentiating and substituting gives the input--output equation
\begin{equation*}
    \ddot{X}_1 + (a+d)\dot{X}_1 + (ad-bc)X_1 = 0.
\end{equation*}
In the case of ideal data, which structural identifiability considers, $X_1$ and $\dot{X}_1$ take all combinations of values independently of each other. The input--output equation must hold for every such combination, which determines its coefficients, $a+d$ and $ad-bc$, uniquely (making them identifiable quantities). These combinations do not uniquely determine $a$, $b$, $c$, and $d$, so the individual parameters are not structurally identifiable from observations of $X_1$ alone. If $d$ were known from other sources, the parameter $a$ and the product $bc$ would become identifiable, but $b$ and $c$ would remain non-identifiable.

Unknown initial conditions, for example in the above example the unobserved variable $X_2$'s initial value $X_2(0)$, can be considered parameters to which structural identifiability can be applied. For brevity, in this work, we will not assess the identifiability of unknown initial conditions.

In this work the unknowns also include the observation map $h$. For instance, if $X_1$ in the example above were measured through $Y=h(X_1)$, then substituting $X_1=g(Y)$ with $g=h^{-1}$ and dividing by $g'(Y)$ gives the input--output equation
\begin{equation*}
    \ddot{Y} + \frac{g''(Y)}{g'(Y)}\dot{Y}^2 + (a+d)\dot{Y} + (ad-bc)\frac{g(Y)}{g'(Y)} = 0.
\end{equation*}
In general, since $h$ has a nowhere-zero derivative it is invertible, so that the observed state variable can be written as $h^{-1}(Y)$ and eliminated from the model equations in the same manner. The coefficients of the resulting input--output equation are then expressions in $h^{-1}$ and its derivatives, and matching them across two choices of the unknowns compares functions rather than numbers. Also note that, as described in Section~\ref{section:obs_map_equation_def}, we assume $h$ to be smooth (which means that the derivatives of $g$, which appear in the input--output equation, are well-defined).


\section{Results} \label{section:results}

Our results are divided into three parts. First, we consider models with a single state variable (Section~\ref{section:1var_models}). Next, we show how differential algebra can determine the identifiability of more general models, using examples from three scientific domains (Section~\ref{section:diff_algebra}). Finally, we consider models in which the observed quantities have dynamics of their own, and so are not read off the state directly (Section~\ref{section:reporters}).

\subsection{Single-variable models} \label{section:1var_models}

We first consider a simple production--decay model. It illustrates several combinations of identifiability and non-identifiability in the core parameters and observation map that can occur. We will next consider a general model of a single variable and use it to prove general results concerning identifiability of single-variable models with unknown observation maps.

\subsubsection{Observations from a production--decay model} \label{section:1var_models_proddeg}

Consider a quantity $X > 0$ that is produced at a constant rate $p\geq0$ and decays at a rate $d>0$ proportional to its abundance. This yields the model
\begin{align}
    \dot{X} &= p - dX,          \label{eq:prod_deg_dynamics} \\
    Y       &= h(X),            \label{eq:prod_deg_obs}
\end{align}
where $Y$ is the observed quantity and $h:\mathbb{R}\to\mathbb{R}$ is an unknown map from the system state to the observable. We assume that $h$ is differentiable with nowhere-zero derivative on the state space. We consider three versions of the model, each of which illustrates a different identifiability pattern. Proofs are provided in Appendix~\ref{appendix:proofs_single_variable_models}.

\begin{proposition} \label{prop:prod_deg_decay}
Consider the model with decay only (i.e. $p\equiv0$), which has dynamics
\begin{equation}
    \dot{X} = -dX. \label{eq:prod_deg_decay_dynamics}
\end{equation}
Then both $d$ and $h$ are structurally non-identifiable. Specifically, given any admissible pair $(h_1,d_1)$, every pair $(h_2,d_2)$ of the form
\begin{align}
    d_2 &\in (0,\infty), \notag \\
    h_2(x) &= h_1\!\left(\left(\frac{x}{k}\right)^{d_1/d_2}\right), \label{eq:prod_deg_decay_nonident}
\end{align}
for any choice of $k>0$, yields identical observed dynamics to $(h_1,d_1)$.
\end{proposition}

\begin{proposition} \label{prop:prod_deg_known_production}
Consider the model with known production rate $p\equiv1$, which has dynamics
\begin{equation}
    \dot{X} = 1-dX. \label{eq:prod_deg_known_production_dynamics}
\end{equation}
Then both $d$ and $h$ are structurally identifiable.
\end{proposition}

\begin{proposition} \label{prop:prod_deg_unknown_production}
Consider the model with an unknown $p>0$, which has dynamics
\begin{equation}
    \dot{X} = p-dX. \label{eq:prod_deg_unknown_production_dynamics}
\end{equation}
Then $d$ is structurally identifiable, whereas $p$ and $h$ are structurally non-identifiable. Specifically, given any admissible choice $(h_1,p_1)$, every choice $(h_2,p_2)$ of the form
\begin{align}
    p_2 &\in (0,\infty), \notag \\
    h_2(x) &= h_1\!\left(\frac{p_1x}{p_2}\right), \label{eq:prod_deg_unknown_production_nonident}
\end{align}
yields identical observed dynamics to $(h_1,p_1)$.
\end{proposition}

These three cases show that the identifiability of the observation map and that of the model parameters need not coincide. With $p\equiv0$, neither can be identified. With $p\equiv1$, both the decay rate and the observation map are identifiable. The most notable case is that in which $p$ is unknown: the decay rate $d$ is identifiable even though the observation map $h$ and the parameter $p$ are not. It is worth noting that just because $h$ is non-identifiable does not mean we cannot learn anything about the model. For example, in the above examples it is identifiable up to some form of rescaling.


\subsubsection{Observations from a general single-variable model} \label{section:1var_models_general}

We next consider the general single-variable model
\begin{align}
    \dot{X} &= F(X;\vec{p}), \label{eq:1var_models_general_dynamics} \\
    Y &= h(X), \label{eq:1var_models_general_obs}
\end{align}
where $\mathcal{D}=(X_{\min},X_{\max})$ is a known interval on which the state $X$ is considered, $F$ is a known functional form (which, for any initial condition, keeps $X(t)\in\mathcal{D}$ for all $t\geq0$), and $\vec{p}$ and $h$ are unknown with $h$ differentiable with nowhere-zero derivative on $\mathcal{D}$. We assume that $\vec{p}$ is structurally identifiable when $X$ is observed directly, i.e. when $h(X) = X$ is known. We give two general results for this model, both of which are proved in Appendix~\ref{appendix:proofs_single_variable_models}.

\begin{proposition} \label{prop:1var_general_nonidentifiability}
Consider the model in Equations~\eqref{eq:1var_models_general_dynamics}--\eqref{eq:1var_models_general_obs}. Suppose that $F(X;\vec{p})\neq0$ for all $X\in\mathcal{D}$ and every admissible $\vec{p}$, and that no solution leaves $\mathcal{D}$ in finite backward time. Then $h$ and $\vec{p}$ are structurally non-identifiable.
\end{proposition}

\begin{proposition} \label{prop:1var_general_identifiability}
Consider the model in Equations~\eqref{eq:1var_models_general_dynamics}--\eqref{eq:1var_models_general_obs}. Suppose that it is known whether $h$ is increasing or decreasing, that every steady state is simple (i.e. $F'(X^*;\vec p)\neq0$ wherever $F(X^*;\vec p)=0$), and that their number and eigenvalues uniquely determine $\vec{p}$. If some solution leaves $\mathcal{D}$ in finite backward time, then both $\vec{p}$ and $h$ are globally structurally identifiable.
\end{proposition}

These propositions can be applied directly to the production--decay model (Propositions~\ref{prop:prod_deg_decay}--\ref{prop:prod_deg_unknown_production}). Proposition~\ref{prop:1var_general_nonidentifiability} applies to the decay-only model: it has no steady state in $(0,\infty)$, and no solution leaves $(0,\infty)$ in finite backward time. The model with $p\equiv1$, by contrast, has a single steady state with eigenvalue $-d$. Since this eigenvalue (trivially) determines the model's single parameter, and because a solution leaves $(0,\infty)$ in finite backward time, Proposition~\ref{prop:1var_general_identifiability} shows that the model is fully identifiable. However, when the production rate is unknown, the eigenvalue $-d$ no longer fully determines the parameter vector $\vec{p}$. Proposition~\ref{prop:1var_general_identifiability} therefore does not apply and cannot be used to determine identifiability.

These results can intuitively be understood by considering that an unknown $h$ acts as an unknown rescaling of the $X$-axis, but leaves the $t$-axis untouched. Hence, data can only reveal model properties that can be inferred from the temporal dynamics. Here, close to any steady state, $h$ is approximately linear. This means that $X$ and $Y$ settle towards their steady state values at the same exponential rates. The eigenvalues can be inferred from these rates, and hence we can infer anything that is determined by the eigenvalues.

\begin{tcolorbox}[colback=black!3,colframe=black!50,boxrule=0.5pt,arc=1mm]
\begin{remark}
Many results throughout this work are only applicable when $h'$ is nowhere zero. In practice, this means that $h$ is strictly monotonic (either increasing or decreasing), which is a natural assumption in many applications. This assumption also helps distinguish the behaviour of the dynamics from that of the observation map. Since $\dot Y=h'(X)\dot X$, allowing $h'$ to vanish means that $Y$ can be stationary even while $X$ is changing; when $h'\neq0$, a stationary observation instead implies that the state itself is stationary. We note that this monotonicity requirement is distinct for unknown functions occurring as observation maps. For example, in \cite{loman_struct_func_2026}, identifiability could be determined for non-monotonic functions occurring within the model dynamics.
\end{remark}
\end{tcolorbox}

\subsection{Determining structural identifiability of observation maps using differential algebra} \label{section:diff_algebra}
Next, we show how differential algebra (Section~\ref{section:differential_algebra}) can be used to determine whether an observation map is identifiable. We consider a linear decay chain, an epidemic model, and a predator--prey model, which between them exhibit a range of different identifiability patterns. Proofs not provided in the main text are given in Appendix~\ref{appendix:proofs_differential_algebra}.

\subsubsection{Linear radioactive decay measured through a faulty detector} \label{section:da_chain}

We first consider a simple model example: a radioactive species $N_1$ decays into a second species $N_2$, which then decays into unobserved products. The abundance of $N_2$ is monitored with a detector carrying an uncharacterised fault, so that the recorded count rate is an unknown function of its true abundance. This gives the model
\begin{align}
    \dot{N}_1 &= -\lambda_1N_1,              \label{eq:chain_dynamics_1} \\
    \dot{N}_2 &= \lambda_1N_1-\lambda_2N_2,  \label{eq:chain_dynamics_2} \\
    Y         &= h(N_2),                     \label{eq:chain_obs}
\end{align}
with $N_1,N_2>0$, distinct unknown rates $\lambda_1,\lambda_2>0$, and $h$ unknown with nowhere-zero derivative.

\begin{proposition} \label{prop:chain}
Consider the model in Equations~\eqref{eq:chain_dynamics_1}--\eqref{eq:chain_obs}. Then $\lambda_1$ and $\lambda_2$ are globally structurally identifiable, whereas $h$ is structurally non-identifiable. Specifically, given any admissible choice $h_1$, the choices giving the same observed outputs are exactly
\begin{equation}
    h_2(x) = h_1(x/k), \qquad k>0. \label{eq:chain_nonident}
\end{equation}
\end{proposition}

\begin{proof}
To generate an input--output equation in the observed variable only, we first extract $N_1$ from Equation~\eqref{eq:chain_dynamics_2}, and also differentiate both sides to gain an expression for $\dot{N}_1$:
\begin{align}
    N_1 &=\frac{1}{\lambda_1}(\dot{N}_2+\lambda_2N_2), \label{eq:chain_N1_eq} \\
    \dot{N}_1 &=\frac{1}{\lambda_1}(\ddot{N}_2+\lambda_2\dot{N}_2). \notag
\end{align}
Inserting these into Equation~\eqref{eq:chain_dynamics_1} gives
\begin{align*}
    \frac{1}{\lambda_1}(\ddot{N}_2+\lambda_2\dot{N}_2) &=-\lambda_1\frac{1}{\lambda_1}(\dot{N}_2+\lambda_2N_2) \Leftrightarrow \\
    \ddot{N}_2 &= -(\lambda_1+\lambda_2)\dot{N}_2-\lambda_1\lambda_2N_2.
\end{align*}
Next, write $g=h^{-1}$, so that $N_2=g(Y)$, and hence $\dot{N}_2=g'(Y)\dot{Y}$ and $\ddot{N}_2=g''(Y)\dot{Y}^2+g'(Y)\ddot{Y}$. Inserting these expressions and dividing by $g'(Y)$ gives
\begin{equation}
    \ddot{Y}=-\frac{g''(Y)}{g'(Y)}\dot{Y}^2-(\lambda_1+\lambda_2)\dot{Y}
    -\lambda_1\lambda_2\frac{g(Y)}{g'(Y)}. \label{eq:chain_proof_io}
\end{equation}
This is the input--output equation. It only depends on the observed variable $Y$ and the coefficients correspond to the identifiable quantities. This gives us that $g''(Y)/g'(Y)$ is identifiable. Assume that two choices $g_1$ and $g_2$ give identical output dynamics. Then
\begin{align*}
    \frac{g_1''(Y)}{g_1'(Y)} &= \frac{g_2''(Y)}{g_2'(Y)} \Leftrightarrow \\
    g_2''g_1'-g_2'g_1'' &= 0 \Leftrightarrow \\
    \left(\frac{g_2'}{g_1'}\right)' &= 0 \Leftrightarrow \\
    g_2 &= kg_1 + b,
\end{align*}
for some constants $k\neq0$ and $b$. In the ideal data case considered by structural identifiability, $N_2$ takes every value in $(0,\infty)$. This means that both $g_1$ and $g_2$ have range $(0,\infty)$, whereas the range of $kg_1+b$ is $(-\infty,b)$ when $k<0$ and $(b,\infty)$ when $k>0$. Hence $k>0$ and $b=0$, so
\begin{equation}
    g_2(Y)=kg_1(Y), \qquad k>0, \label{eq:chain_proof_g}
\end{equation}
which is equivalent to $h_2(x)=h_1(x/k)$.

Since $g_2=kg_1$ we also have $g_2'=kg_1'$, and hence
\begin{equation*}
    \frac{g_2(Y)}{g_2'(Y)}=\frac{kg_1(Y)}{kg_1'(Y)}=\frac{g_1(Y)}{g_1'(Y)} ,
\end{equation*}
so that $g(Y)/g'(Y)$ is identifiable, even though $g$ itself is not. The constant coefficient of Equation~\eqref{eq:chain_proof_io} is $-\lambda_1\lambda_2\,g(Y)/g'(Y)$, and dividing it by $g(Y)/g'(Y)$ shows that $\lambda_1\lambda_2$ is identifiable. The coefficient of $\dot{Y}$ is $-(\lambda_1+\lambda_2)$, so $\lambda_1+\lambda_2$ is identifiable as well. We hence know the sum and the product of the two rates, and so the pair $\{\lambda_1,\lambda_2\}$, but not which of the two is $\lambda_2$. One further identifiable quantity is required to prove full identifiability.

We obtain it from Equation~\eqref{eq:chain_N1_eq}, which we have so far used only as an identity, and not as a constraint. It holds with $N_1>0$. Substituting $N_2=g(Y)$ and $\dot{N}_2=g'(Y)\dot{Y}$ into the inequality gives
\begin{align*}
    N_1 &> 0 \Leftrightarrow \\
    \dot{N}_2+\lambda_2N_2 &> 0 \Leftrightarrow \\
    g'(Y)\dot{Y}+\lambda_2g(Y) &> 0.
\end{align*}
At each $Y$, the observed values of $\dot{Y}$ are therefore exactly those lying on one side of the curve
\begin{equation}
    \dot{Y}=-\lambda_2\frac{g(Y)}{g'(Y)}. \label{eq:chain_proof_endpoint}
\end{equation}
They lie above this curve when $g'>0$, and below it when $g'<0$. The curve is therefore the boundary of the observed values of $\dot{Y}$, so it is a feature of the observed dynamics, making the boundary $\lambda_2\,g(Y)/g'(Y)$ identifiable. Since $g(Y)/g'(Y)$ is identifiable and non-zero, we have that $\lambda_2$ is identifiable. Identifiability of $\lambda_1$ then follows from the identifiability of $\lambda_1+\lambda_2$. Both rates are hence globally structurally identifiable.

Conversely, for every $k>0$, scaling both $N_1$ and $N_2$ by $k$ preserves the dynamics and is absorbed by $h_2(x)=h_1(x/k)$. This gives exactly the family in Equation~\eqref{eq:chain_nonident}.
\end{proof}

Here we see the same phenomenon as in the production--decay model with an unknown production rate: the core model parameters are all identifiable, even though the observation map is not. This is useful, as the core system is typically the one of interest. Our results show that we are able to determine this core system fully, even when the device with which we measure the system is subject to systematic, but unknown, faults (which prevent us from knowing exactly how the system state relates to measurements).

We note that for this model, $h$ is identifiable up to an unknown rescaling. However, let us assume that we have access to a calibration sample, i.e. that we can learn $h(x_0)$ for some specific value of $x_0$. In the above example this could be a sample where a known quantity of $N_2$ is measured to yield a known value of $h(N_2)$. Under these circumstances, we have full identifiability (proof in Appendix~\ref{appendix:proof_chain_calibration}). This is a relatively common case in many practical experimental workflows, where such calibration samples are available and used to calibrate various machines. This demonstrates how non-dynamic data can be an important source of information for improving identifiability of dynamic models.

\subsubsection{An SIR model where infected individuals are measured through an unknown function} \label{section:da_sir}

Next, we consider the spread of an infectious disease, as described by the classic SIR model~\cite{anderson_infectious_1991}. Unlike in the standard model, however, we can only measure the infected cases indirectly, for example through wastewater measurements~\cite{peccia_measurement_2020,huisman_wastewater-based_2022,schmid_wastewater-informed_2026}. This gives the model
\begin{align}
    \dot{S} &= -\beta SI,           \label{eq:sir_dynamics_s} \\
    \dot{I} &= \beta SI - \gamma I, \label{eq:sir_dynamics_i} \\
    \dot{R} &= \gamma I,            \label{eq:sir_dynamics_r} \\
    Y       &= h(I),                \label{eq:sir_obs}
\end{align}
with $S > 0$, $I > 0$, $R \geq 0$, and the unknown parameters are $\beta>0$ and $\gamma>0$. We assume that the unknown $h$ is differentiable with $h' \neq 0 $ everywhere. 

\begin{proposition} \label{prop:sir}
Consider the model in Equations~\eqref{eq:sir_dynamics_s}--\eqref{eq:sir_obs}. Here, $\gamma$ is globally structurally identifiable, whereas $\beta$ and $h$ are structurally non-identifiable. Specifically, given any admissible choice $(h_1,\beta_1)$, the choices giving the same observed outputs are exactly
\begin{align}
    \beta_2  &= \beta_1/c, \notag \\
    h_2(x)   &= h_1(x/c), \qquad c>0. \label{eq:sir_nonident}
\end{align}
If we also know the total population size, $N=S+I+R$ (which remains constant), then $\beta$, $\gamma$ and $h$ are all globally structurally identifiable.
\end{proposition}

This shows that we can only fully learn the model when we track the relative abundance of the infection throughout the population, rather than the absolute numbers. As the total population size is often known, however, this will rarely be an actual constraint. Generally, this implies that it is possible to learn at least simple epidemiological models from indirect measurements alone.

\begin{figure}[!htb]
    \centering
    \includegraphics[width=0.98\linewidth]{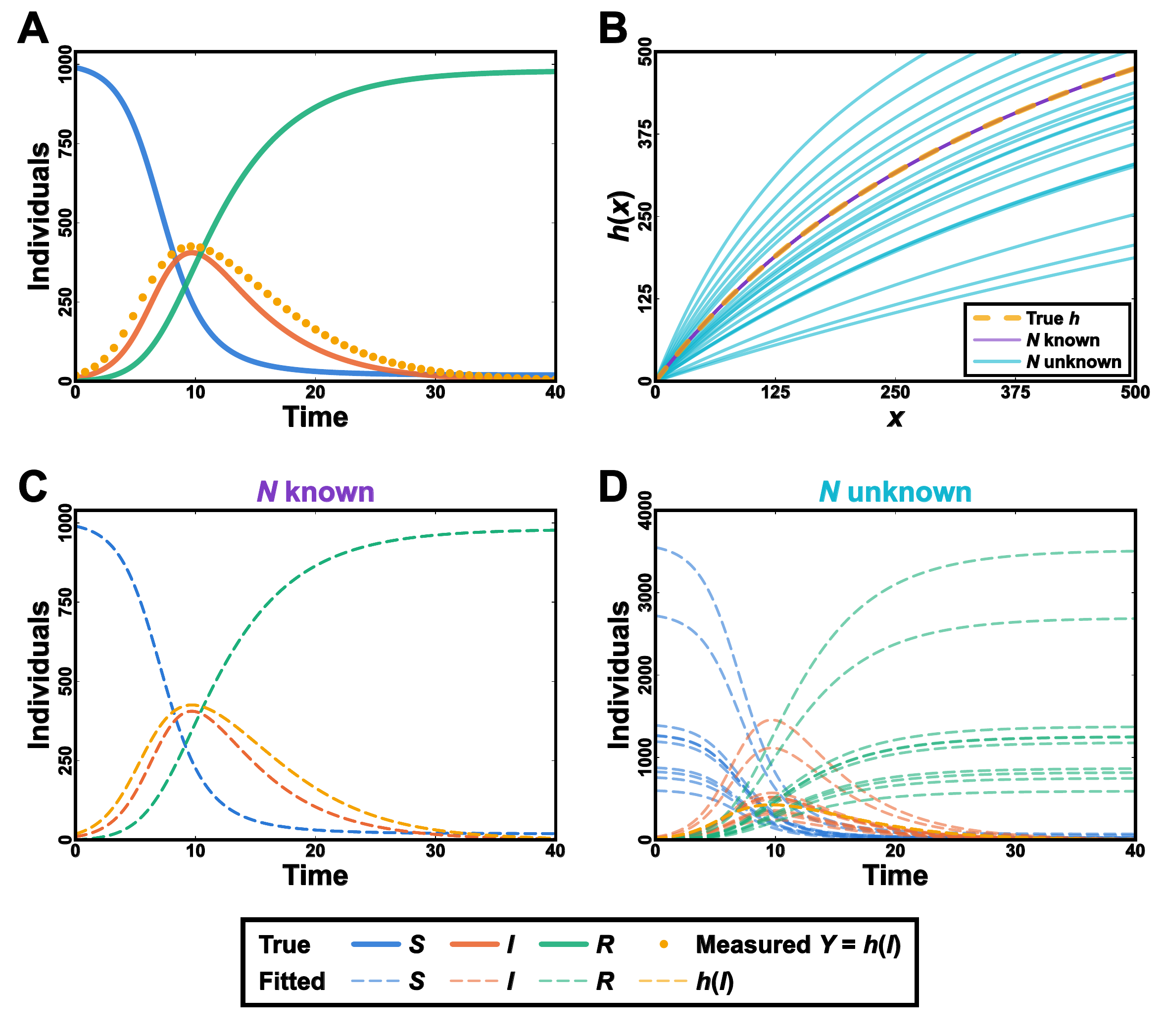}
    \caption{\textbf{Knowing the total population size makes the SIR model with an unknown observation map identifiable.} (A) Ground-truth trajectories of $S$, $I$, and $R$ (solid lines), and the synthetic measurements of $Y=h(I)$ (dots) to which the model is fitted. (B) The ground-truth observation map $h$ (orange, dashed), and the maps recovered by the 20 best fits when the total population size $N$ is known (purple) and unknown (cyan). Only where we know the total population size is the observation map identified correctly (when we do not, any rescaled version can reproduce the dynamics). (C) From repeated randomised fits of the model where $N$ is known, we extract the 10 best fits and plot the trajectories of $S$, $I$, $R$, and $h(I)$. All fits coincide with the true trajectories. (D) Repeat of C, but where $N$ is unknown. Here, while all $h(I)$ trajectories align with the measurements, there is great variability for the $S$, $I$, and $R$ trajectories.  Each trajectory retains the shape of the original dynamics, but with an arbitrary rescaling (as predicted by Proposition~\ref{prop:sir}). The data generation and fitting procedures are described in Appendix~\ref{appendix:numerical_implementation}.}
    \label{fig:sir_fits}
\end{figure}

In Figure~\ref{fig:sir_fits} we demonstrate the above result numerically. From a synthetic data set (generated as described in Appendix~\ref{appendix:sir_data_generation}) we fit $\beta$, $\gamma$, and $h$ for the SIR model described in Equations~\eqref{eq:sir_dynamics_s}--\eqref{eq:sir_obs} (fitting workflow described in Appendix~\ref{appendix:sir_fitting}). We consider two cases, one where the total population size is known, and one where it is not. In the first case, we can recover the true dynamics, as well as $h$. In the second case, the dynamics and $h$ can be recovered only up to an unknown scaling factor.

\subsubsection{A Lotka--Volterra model where predators are measured through an unknown function} \label{section:da_lotka_volterra}

Finally, we consider a Lotka--Volterra predator--prey model~\cite{alfred_jlotka_elements_1925,volterra_fluctuations_1926} in which the predator's abundance is measured through an unknown function. This is a common situation in ecology, where populations sometimes cannot be observed directly, and are instead monitored through proxies such as collected faeces~\cite{marques_estimating_2001}, tracks, or indirect observations like recording mating calls~\cite{marques_estimating_2013}. Specifically, we consider the following model:
\begin{align}
    \dot{N} &= rN - aNP, \label{eq:lotka_volterra_dynamics_prey} \\
    \dot{P} &= NP - mP,  \label{eq:lotka_volterra_dynamics_pred} \\
    Y       &= h(P),      \label{eq:lotka_volterra_obs}
\end{align}
with prey and predator abundances $N>0$ and $P>0$, unknown parameters $r>0$, $a>0$ and $m>0$, and the unknown observation map $h$ assumed to have nowhere-zero derivative.

\begin{proposition} \label{prop:lotka_volterra}
Consider the model in Equations~\eqref{eq:lotka_volterra_dynamics_prey}--\eqref{eq:lotka_volterra_obs}. Then $r$ and $m$ are globally structurally identifiable, whereas $a$ and $h$ are structurally non-identifiable. Specifically, given any admissible choice $(h_1,a_1)$, the choices giving the same observed outputs are exactly
\begin{align}
    a_2    &= a_1/c, \notag \\
    h_2(x) &= h_1(x/c), \qquad c>0. \label{eq:lotka_volterra_nonident}
\end{align}
\end{proposition}

Here, indirect observations are enough to learn some information about the core system, but not enough to determine its full parameterisation. This non-identifiability holds even if we, in addition to measuring the predator population through one (or even several) unknown map(s), can measure the prey population directly (Appendix~\ref{appendix:lotka_volterra_multiple_maps}). This distinguishes the Lotka--Volterra model from the SIR model, where full structural identifiability was possible. Here, we instead have a fundamental limitation on what can be learnt through unknown observation maps.

\subsection{Observable reporters produced at an unknown rate} \label{section:reporters}

In the sections above, we have assumed that, although the observed quantity is seen through an unknown observation map, the observations themselves are instantaneous in time. This assumption does not always hold. Consider, for example, the observation of an infectious disease outbreak through wastewater markers. The model in Equations~\eqref{eq:sir_dynamics_s}--\eqref{eq:sir_obs} assumes that the dynamics of the outbreak are much slower than those of the wastewater markers. This is not always expected to be true. In this section, we instead consider cases where the observed quantity is a reporter with dynamics of its own. Specifically, we consider the following model:
\begin{align}
    \dot{\vec{X}} &= \vec{F}(\vec{X};\vec{p}),                              \label{eq:reporter_ode} \\
    \dot{R}_j     &= f_j(\vec{u}_j) - d_jR_j, \qquad j = 1,\ldots,k,        \label{eq:reporter_marker} \\
    Y_j           &= R_j,                                                   \label{eq:reporter_obs}
\end{align}
where $f_j$ is an unknown production rate acting on a sub-vector $\vec{u}_j$ of the state, and $d_j>0$ is the reporter's decay rate. The unknowns are $\vec{f}=(f_1,\ldots,f_k),\vec{p} = (p_1, \ldots, p_l), \vec{d} = (d_1, \ldots, d_k)$. A well-known class of systems following these dynamics are genetic circuits observed through fluorescent reporter proteins. Here, the production rate of the reporter is an unknown function of the activating components, while the reporter itself can be assumed to decay at a linear rate~\cite{ronen_assigning_2002,de_jong_experimental_2010}. 

Strictly, the model in Equations~\eqref{eq:reporter_ode}--\eqref{eq:reporter_obs} is not an observation map in the sense of Equations~\eqref{eq:base_model_ode}--\eqref{eq:base_model_obs}, but rather a standard functional inference problem of the type studied in~\cite{loman_struct_func_2026}. However, besides being phenomenologically equivalent to what we study here, it can also be reduced to the same mathematical problem, as we show in the following two propositions (proved in Appendix~\ref{appendix:proofs_reporters}).

\begin{proposition} \label{prop:reporter_reduction}
Suppose the decay rates $\vec{d}$ are known. Then the unknown function $f_j$ is structurally identifiable if and only if it is structurally identifiable in the corresponding observation-map model
\begin{align}
    \dot{\vec{X}} &= \vec{F}(\vec{X};\vec{p}), \notag \\
    Y_j &= f_j(\vec{u}_j), \qquad j = 1,\ldots,k. \label{eq:reporter_reduced}
\end{align}
\end{proposition}

The result follows directly from integrating Equation~\eqref{eq:reporter_marker} (Appendix~\ref{appendix:proof_reporter_reduction}). That the decay rates are known is not necessarily a strong assumption. For example, in many biological systems, decay is dominated by dilution due to cell division~\cite{alon_introduction_2019}, which is often either known or possible to measure directly. The following proposition also gives a condition under which $d_j$, and hence $f_j$, are structurally identifiable.

\begin{proposition} \label{prop:reporter_decay}
Consider the model in Equations~\eqref{eq:reporter_ode}--\eqref{eq:reporter_obs}. Suppose that, for any initial value of $\vec{X}$, all non-negative reporter values are admissible initial conditions of the system. If each unknown production rate $f_j$ is known to be non-negative, then each decay rate $d_j$ is globally structurally identifiable. This also implies that Proposition~\ref{prop:reporter_reduction} then applies.
\end{proposition}

The two conditions should be widely applicable. I.e. non-negativity of $\vec{f}$ is expected across many cases. Simultaneously, for most systems, there is little stopping any pairing of initial conditions $\vec{X}(0)$ and $\vec{R}(0)$. For example, in the case of fluorescent reporters, this means that for any potential initial state of the system $\vec{X}(0)$, any potential initial abundance of the reporter $\vec{R}(0)$ is admissible, which is a realistic assumption. 

Together, these two propositions show how identifiability for models of the form in Equations~\eqref{eq:reporter_ode}--\eqref{eq:reporter_obs} can be reduced to the observation map problem studied throughout this work. Next, we consider how this can be done in practice, focusing on the case of fluorescent reporters, which are widely used in biology.

\subsubsection{A reporter of a decaying protein} \label{section:reporters_decay}

First, we consider the simplest case: a protein that decays at a linear rate and is observed through a fluorescent reporter. This gives the model
\begin{align}
    \dot{X} &= -aX,       \label{eq:reporter_decay_ode} \\
    \dot{R} &= f(X) - dR, \label{eq:reporter_decay_marker} \\
    Y       &= R,          \label{eq:reporter_decay_obs}
\end{align}
where $X>0$ is the protein and $a>0$ its degradation rate, $R$ is the reporter and $d>0$ its decay rate, and $f$ is the unknown production rate, taken to be positive and differentiable with nowhere-zero derivative. The initial value of $X$ is unknown.

\begin{proposition} \label{prop:reporter_decay_signal}
Consider the model in Equations~\eqref{eq:reporter_decay_ode}--\eqref{eq:reporter_decay_obs}. Then $d$ is globally structurally identifiable, whereas $a$ and $f$ are structurally non-identifiable. Specifically, given any admissible choice $(f_1,a_1)$, the choices giving the same observed outputs are exactly
\begin{align}
    a_2    &\in(0,\infty), \notag \\
    f_2(x) &= f_1\bigl((x/c)^{a_1/a_2}\bigr), \qquad c>0. \label{eq:reporter_decay_nonident}
\end{align}
\end{proposition}
The proof (Appendix~\ref{appendix:proof_reporter_decay_signal}) combines Propositions~\ref{prop:reporter_reduction} and~\ref{prop:reporter_decay} to reduce this model to the decay-only model of Section~\ref{section:1var_models}, from which the non-identifiability follows.

\subsubsection{A reporter of a genetic oscillator} \label{section:reporters_goodwin}
Next, we consider a more complicated model: a Goodwin oscillator~\cite{goodwin_oscillatory_1965,gonze_goodwin_2013}. This is a well-known model of a biological oscillator driven by a delayed negative feedback loop. Here, we study whether the Goodwin oscillator model is identifiable when it can be measured through a fluorescent reporter only. This gives the following model:
\begin{align}
    \dot{X}_1 &= a/(1+X_3^m) - b_1X_1, \label{eq:reporter_goodwin_x1} \\
    \dot{X}_2 &= X_1 - b_2X_2,         \label{eq:reporter_goodwin_x2} \\
    \dot{X}_3 &= X_2 - b_3X_3,         \label{eq:reporter_goodwin_x3} \\
    \dot{R}   &= f(X_3) - dR,          \label{eq:reporter_goodwin_marker} \\
    Y         &= R,                     \label{eq:reporter_goodwin_obs}
\end{align}
where $X_1$, $X_2$ and $X_3$ are the components of the oscillator, none of them observed and all positive, where $a$, $m$, $b_1$, $b_2$, $b_3$ and $d$ are unknown positive parameters, and where $f$ is an unknown production rate, positive and differentiable with nowhere-zero derivative. The following proposition (proved in Appendix~\ref{appendix:proof_reporter_goodwin}) shows that the model is fully identifiable, even though we can only observe the reporter.

\begin{proposition} \label{prop:reporter_goodwin}
Consider the model in Equations~\eqref{eq:reporter_goodwin_x1}--\eqref{eq:reporter_goodwin_obs}. Then all of the unknowns, that is $a$, $m$, $b_1$, $b_2$, $b_3$, $d$ and the production rate $f$, are globally structurally identifiable.
\end{proposition}

\begin{tcolorbox}[colback=black!3,colframe=black!50,boxrule=0.5pt,arc=1mm]
\begin{remark}
In this section, we have shown how some directly observed models with unknown functions in the dynamics can be turned into models with unknown observation maps, whose structural identifiability can then be determined. The reverse step is also possible: any model of the form in Equations~\eqref{eq:base_model_ode}--\eqref{eq:base_model_obs} can be written as a model with unknown functions in the dynamics. Differentiating both sides of Equation~\eqref{eq:base_model_obs} gives
\begin{align*}
    \dot{\vec{X}} &= \vec{F}(\vec{X};\vec{p}), \\
    \dot{Y}_j     &= \nabla h_j(\vec{X})\cdot\vec{F}(\vec{X};\vec{p}), \qquad j=1,\ldots,k,
\end{align*}
with initial values linked through $\vec{Y}(0)=\vec{h}(\vec{X}(0))$. For a single state variable and a single observation map, this reads $\dot{Y}=h'(X)F(X;\vec{p})$. The structural identifiability of $\vec{p}$ and $\vec{h}$ can then also be determined from this form (using methodology described in~\cite{loman_struct_func_2026}). 
\end{remark}
\end{tcolorbox}


\section{Discussion} \label{section:discussion}

In this work, we considered differential equation models subject to observations where the map from the system state to observations is unknown. We asked when these models exhibit identifiability, both considering the core dynamics and the observation mapping. Specifically, we both developed general results on when these models exhibit identifiability, and demonstrated a general methodology for determining identifiability for individual models. This is based on the differential algebra approach, an established method for determining parameter identifiability~\cite{ljung_global_1994,audoly_global_2001,bellu_daisy_2007}, which we recently also introduced as an approach for determining the identifiability of unknown functions~\cite{loman_struct_func_2026}. 

We specifically considered multiple relevant examples in which such unknown maps occur naturally. These include indirect measurements of systems using detectors subject to uncharacterised faults, measurements of infectious disease outbreaks (for example through markers in wastewater), interacting animal populations (through, for example, faeces collection or records of bird song), and measurements of fluorescent markers. All of these are real and highly relevant situations in which such unknown observation maps arise. Notable results include that SIR models observed through such maps are structurally identifiable, suggesting that infectious disease outbreaks can be modelled from this kind of data alone.

We also considered models observed through reporters, i.e.~observables with dynamics of their own that do not feed back into the core system. For reporters produced at an unknown rate and decaying linearly, we showed that these models can often be linked directly to observation-map models. Related models, for instance ones in which the reporter is produced at a linear rate but decays according to an unknown function, are not treated here and are a subject for further study.

Similarly, our work can be generalised in many directions. We have considered a few model examples only. The goal has been to demonstrate how the differential algebra approach can be applied and to show promising example cases for applying unknown observation map models. We have also provided general results (like Propositions~\ref{prop:1var_general_nonidentifiability} and~\ref{prop:1var_general_identifiability}) to show that structural identifiability, or lack thereof, can be determined for general categories of models. Likely, we can obtain much stronger general results on how properties of e.g. the system dynamics $F$ determine when identifiability occurs. Such results are intrinsically useful in themselves, but would also provide additional general insight into the nature of observation maps (for example, when are non-identifiabilities relatively simple rescalings of $h$, and do more complicated non-identifiabilities occur). Again, however, our goal has not been to provide an all-encompassing theory of when such observation maps are identifiable, but rather to demonstrate the kinds of results one can obtain, and ideas that can take one there.

Similarly, the fact that a map is structurally identifiable does not mean that it is practically identifiable. Work in~\cite{loman_functional_2025} has shown how to assess such practical functional identifiability. By applying these methods to models with unknown observation maps, one could survey under what settings these models exhibit practical identifiability. Here, our work forms a solid basis for showing that these kinds of map can, in theory, be learnt, and that they are a worthwhile avenue for future research.

The current rise of machine learning within modelling has increasingly promoted the integration of data-driven modelling techniques with classic mechanistic mathematical modelling~\cite{noordijk_rise_2024,karniadakis_physics-informed_2021,brunton_discovering_2016,chen_neural_2018}. Accounting for unknown observation maps is one natural such extension. Until recently, such maps could not be represented as unknown functions, so how the observations depend on the system had to be fully specified in advance. A wrongly specified map would then introduce model misspecification errors. While hybrid models offer a principled approach for modelling such systems, the approach has received little attention in this context. By demonstrating how to assess identifiability for this class of models, and showing that they often permit either full or extensive identification of system dynamics, we hope to further promote their application to real-world systems.




\section*{Acknowledgements} \label{section:acknowledgements}
TEL and REB are supported by a grant from the Simons Foundation (MP-SIP-00001828). For the purpose of open access, the author has applied a CC BY public copyright licence to any author accepted manuscript arising from this submission. AI tools were used to assist with code development, manuscript editing, and refinement of the exposition and mathematical arguments. All final text and code have been verified, understood, and edited by the authors before publication.


\bibliographystyle{plainnat}
\bibliography{references}

@book{alfred_jlotka_elements_1925,
	title = {Elements of {Physical} {Biology}},
	publisher = {Williams and Wilkins Company},
	author = {Lotka, Alfred J.},
	year = {1925},
}

@book{alon_introduction_2019,
	title = {An {Introduction} to {Systems} {Biology}: {Design} {Principles} of {Biological} {Circuits}},
	edition = {2nd},
	doi = {10.1201/9780429283321},
	publisher = {Chapman and Hall/CRC},
	author = {Alon, Uri},
	year = {2019},
}

@book{anderson_infectious_1991,
	title = {Infectious {Diseases} of {Humans}: {Dynamics} and {Control}},
	doi = {10.1093/oso/9780198545996.001.0001},
	publisher = {Oxford University Press},
	author = {Anderson, Roy M. and May, Robert M.},
	year = {1991},
}

@article{audoly_global_2001,
	title = {Global identifiability of nonlinear models of biological systems},
	volume = {48},
	doi = {10.1109/10.900248},
	number = {1},
	journal = {IEEE Transactions on Biomedical Engineering},
	author = {Audoly, S. and Bellu, G. and D'Angio, L. and Saccomani, M.P. and Cobelli, C.},
	year = {2001},
	pages = {55--65},
}

@article{bellman_structural_1970,
	title = {On structural identifiability},
	volume = {7},
	doi = {10.1016/0025-5564(70)90132-X},
	number = {3},
	journal = {Mathematical Biosciences},
	author = {Bellman, R. and Åström, K. J.},
	year = {1970},
	pages = {329--339},
}

@article{bellu_daisy_2007,
	title = {{DAISY}: {A} new software tool to test global identifiability of biological and physiological systems},
	volume = {88},
	doi = {10.1016/j.cmpb.2007.07.002},
	number = {1},
	journal = {Computer Methods and Programs in Biomedicine},
	author = {Bellu, Giuseppina and Saccomani, Maria Pia and Audoly, Stefania and D’Angiò, Leontina},
	year = {2007},
	pages = {52--61},
}

@article{bezanson_julia_2017,
	title = {{Julia}: {A} fresh approach to numerical computing},
	volume = {59},
	doi = {10.1137/141000671},
	number = {1},
	journal = {SIAM Review},
	author = {Bezanson, Jeff and Edelman, Alan and Karpinski, Stefan and Shah, Viral B.},
	year = {2017},
	pages = {65--98},
}

@article{brunton_discovering_2016,
	title = {Discovering governing equations from data by sparse identification of nonlinear dynamical systems},
	volume = {113},
	doi = {10.1073/pnas.1517384113},
	number = {15},
	journal = {Proceedings of the National Academy of Sciences},
	author = {Brunton, Steven L. and Proctor, Joshua L. and Kutz, J. Nathan},
	year = {2016},
	pages = {3932--3937},
}

@inproceedings{chen_neural_2018,
	title = {Neural ordinary differential equations},
	booktitle = {Advances in {Neural} {Information} {Processing} {Systems} ({NeurIPS})},
	volume = {31},
	url = {https://arxiv.org/abs/1806.07366},
	author = {Chen, Ricky T. Q. and Rubanova, Yulia and Bettencourt, Jesse and Duvenaud, David},
	year = {2018},
}

@article{cobelli_parameter_1980,
	title = {Parameter and structural identifiability concepts and ambiguities: a critical review and analysis},
	volume = {239},
	doi = {10.1152/ajpregu.1980.239.1.R7},
	number = {1},
	journal = {The American Journal of Physiology},
	author = {Cobelli, C. and DiStefano, J. J.},
	year = {1980},
	pages = {R7--24},
}

@article{de_jong_experimental_2010,
	title = {Experimental and computational validation of models of fluorescent and luminescent reporter genes in bacteria},
	volume = {4},
	doi = {10.1186/1752-0509-4-55},
	journal = {BMC Systems Biology},
	author = {de Jong, Hidde and Ranquet, Caroline and Ropers, Delphine and Pinel, Corinne and Geiselmann, Johannes},
	year = {2010},
	pages = {55},
}

@article{gonze_goodwin_2013,
	title = {The {Goodwin} model: behind the {Hill} function},
	volume = {8},
	doi = {10.1371/journal.pone.0069573},
	number = {8},
	journal = {PLoS ONE},
	author = {Gonze, Didier and Abou-Jaoudé, Wassim},
	year = {2013},
	pages = {e69573},
}

@article{goodwin_oscillatory_1965,
	title = {Oscillatory behavior in enzymatic control processes},
	volume = {3},
	doi = {10.1016/0065-2571(65)90067-1},
	journal = {Advances in Enzyme Regulation},
	author = {Goodwin, Brian C.},
	year = {1965},
	pages = {425--437},
}

@article{huisman_wastewater-based_2022,
	title = {Wastewater-based estimation of the effective reproductive number of {SARS-CoV-2}},
	volume = {130},
	doi = {10.1289/EHP10050},
	number = {5},
	journal = {Environmental Health Perspectives},
	author = {Huisman, Jana S. and Scire, Jérémie and Caduff, Lea and Fernandez-Cassi, Xavier and Ganesanandamoorthy, Pravin and Kull, Anina and Scheidegger, Andreas and Stachler, Elyse and Boehm, Alexandria B. and Hughes, Bridgette and Knudson, Alisha and Topol, Aaron and Wigginton, Krista R. and Wolfe, Marlene K. and Kohn, Tamar and Ort, Christoph and Stadler, Tanja and Julian, Timothy R.},
	year = {2022},
	pages = {057011},
}

@article{karniadakis_physics-informed_2021,
	title = {Physics-informed machine learning},
	volume = {3},
	doi = {10.1038/s42254-021-00314-5},
	number = {6},
	journal = {Nature Reviews Physics},
	author = {Karniadakis, George Em and Kevrekidis, Ioannis G. and Lu, Lu and Perdikaris, Paris and Wang, Sifan and Yang, Liu},
	year = {2021},
	pages = {422--440},
}

@inproceedings{kingma_adam_2015,
	title = {Adam: {A} method for stochastic optimization},
	booktitle = {International Conference on Learning Representations ({ICLR})},
	url = {https://arxiv.org/abs/1412.6980},
	author = {Kingma, Diederik P. and Ba, Jimmy},
	year = {2015},
}

@article{liu_limited_1989,
	title = {On the limited memory {BFGS} method for large scale optimization},
	volume = {45},
	doi = {10.1007/BF01589116},
	number = {1--3},
	journal = {Mathematical Programming},
	author = {Liu, Dong C. and Nocedal, Jorge},
	year = {1989},
	pages = {503--528},
}

@article{ljung_global_1994,
	title = {On global identifiability for arbitrary model parametrizations},
	volume = {30},
	doi = {10.1016/0005-1098(94)90029-9},
	number = {2},
	journal = {Automatica},
	author = {Ljung, Lennart and Glad, Torkel},
	year = {1994},
	pages = {265--276},
}

@misc{loman_functional_2025,
	title = {Functional and parametric identifiability for universal differential equations applied to chemical reaction networks},
	howpublished = {arXiv:2510.14140},
	author = {Loman, Torkel E. and Baker, Ruth E.},
	year = {2025},
}

@misc{loman_struct_func_2026,
	title = {Structural functional identifiability and model discovery in differential equation models},
	author = {Loman, Torkel E. and Browning, Alexander P. and Baker, Ruth E.},
	year = {2026},
	howpublished = {arXiv:2606.30289},
}

@article{marques_estimating_2001,
	title = {Estimating deer abundance from line transect surveys of dung: sika deer in southern {Scotland}},
	volume = {38},
	doi = {10.1046/j.1365-2664.2001.00584.x},
	number = {2},
	journal = {Journal of Applied Ecology},
	author = {Marques, Fernanda F. C. and Buckland, Stephen T. and Goffin, David and Dixon, Camilla E. and Borchers, David L. and Mayle, Brenda A. and Peace, Andrew J.},
	year = {2001},
	pages = {349--363},
}

@article{marques_estimating_2013,
	title = {Estimating animal population density using passive acoustics},
	volume = {88},
	doi = {10.1111/brv.12001},
	number = {2},
	journal = {Biological Reviews},
	author = {Marques, Tiago A. and Thomas, Len and Martin, Stephen W. and Mellinger, David K. and Ward, Jessica A. and Moretti, David J. and Harris, Danielle and Tyack, Peter L.},
	year = {2013},
	pages = {287--309},
}

@article{mogensen_optim_2018,
	title = {Optim: {A} mathematical optimization package for {Julia}},
	volume = {3},
	doi = {10.21105/joss.00615},
	number = {24},
	journal = {Journal of Open Source Software},
	author = {Mogensen, Patrick K. and Riseth, Asbjørn N.},
	year = {2018},
	pages = {615},
}

@article{noordijk_rise_2024,
	title = {The rise of scientific machine learning: a perspective on combining mechanistic modelling with machine learning for systems biology},
	volume = {4},
	doi = {10.3389/fsysb.2024.1407994},
	journal = {Frontiers in Systems Biology},
	author = {Noordijk, Ben and Garcia Gomez, Monica L. and ten Tusscher, Kirsten H. W. J. and de Ridder, Dick and van Dijk, Aalt D. J. and Smith, Robert W.},
	year = {2024},
}

@article{peccia_measurement_2020,
	title = {Measurement of {SARS-CoV-2} {RNA} in wastewater tracks community infection dynamics},
	volume = {38},
	doi = {10.1038/s41587-020-0684-z},
	number = {10},
	journal = {Nature Biotechnology},
	author = {Peccia, Jordan and Zulli, Alessandro and Brackney, Doug E. and Grubaugh, Nathan D. and Kaplan, Edward H. and Casanovas-Massana, Arnau and Ko, Albert I. and Malik, Amyn A. and Wang, Dennis and Wang, Mike and Warren, Joshua L. and Weinberger, Daniel M. and Arnold, Wyatt and Omer, Saad B.},
	year = {2020},
	pages = {1164--1167},
}

@misc{persson_petab_sciml_2026,
	title = {{PEtab} {SciML}: an exchange format for specifying and training dynamic scientific machine learning models},
	author = {Persson, Sebastian and Snelling, Branwen and Philipps, Maren and Weindl, Daniel and Cvijovic, Marija and Hasenauer, Jan and Pathirana, Dilan and Fröhlich, Fabian},
	year = {2026},
	howpublished = {arXiv:2608.20184},
}

@misc{rackauckas_universal_2021,
	title = {Universal {Differential} {Equations} for {Scientific} {Machine} {Learning}},
	howpublished = {arXiv:2001.04385},
	author = {Rackauckas, Christopher and Ma, Yingbo and Martensen, Julius and Warner, Collin and Zubov, Kirill and Supekar, Rohit and Skinner, Dominic and Ramadhan, Ali and Edelman, Alan},
	year = {2021},
}

@article{raue_structural_2009,
	title = {Structural and practical identifiability analysis of partially observed dynamical models by exploiting the profile likelihood},
	volume = {25},
	doi = {10.1093/bioinformatics/btp358},
	number = {15},
	journal = {Bioinformatics},
	author = {Raue, A. and Kreutz, C. and Maiwald, T. and Bachmann, J. and Schilling, M. and Klingmüller, U. and Timmer, J.},
	year = {2009},
	pages = {1923--1929},
}

@article{ronen_assigning_2002,
	title = {Assigning numbers to the arrows: parameterizing a gene regulation network by using accurate expression kinetics},
	volume = {99},
	doi = {10.1073/pnas.152046799},
	number = {16},
	journal = {Proceedings of the National Academy of Sciences},
	author = {Ronen, Michal and Rosenberg, Revital and Shraiman, Boris I. and Alon, Uri},
	year = {2002},
	pages = {10555--10560},
}

@misc{revels_forward-mode_2016,
	title = {Forward-mode automatic differentiation in {Julia}},
	howpublished = {arXiv:1607.07892},
	author = {Revels, Jarrett and Lubin, Miles and Papamarkou, Theodore},
	year = {2016},
}

@article{simpson_parameter_2026,
	title = {Parameter {Identifiability}, {Parameter} {Estimation}, and {Model} {Prediction} for {Differential} {Equation} {Models}},
	volume = {68},
	doi = {10.1137/24M1667968},
	journal = {SIAM Review},
	author = {Simpson, Matthew and Baker, Ruth},
	year = {2026},
}

@article{simpson_profile-wise_2023,
	title = {Profile-{Wise} {Analysis}: {A} profile likelihood-based workflow for identifiability analysis, estimation, and prediction with mechanistic mathematical models},
	volume = {19},
	doi = {10.1371/journal.pcbi.1011515},
	number = {9},
	journal = {PLOS Computational Biology},
	author = {Simpson, Matthew J. and Maclaren, Oliver J.},
	year = {2023},
	pages = {e1011515},
}

@article{volterra_fluctuations_1926,
	title = {Fluctuations in the {Abundance} of a {Species} considered {Mathematically}},
	volume = {118},
	doi = {10.1038/118558a0},
	number = {2972},
	journal = {Nature},
	author = {Volterra, Vito},
	year = {1926},
	pages = {558--560},
}

@article{wieland_structural_2021,
	title = {On structural and practical identifiability},
	volume = {25},
	doi = {10.1016/j.coisb.2021.03.005},
	journal = {Current Opinion in Systems Biology},
	author = {Wieland, Franz-Georg and Hauber, Adrian L. and Rosenblatt, Marcus and Tönsing, Christian and Timmer, Jens},
	year = {2021},
	pages = {60--69},
}

@article{bonnaffe_neural_2021,
	title = {Neural ordinary differential equations for ecological and evolutionary time-series analysis},
	volume = {12},
	doi = {10.1111/2041-210X.13606},
	number = {7},
	journal = {Methods in Ecology and Evolution},
	author = {Bonnaff{\'e}, Willem and Sheldon, Ben C. and Coulson, Tim},
	year = {2021},
	pages = {1301--1315},
}

@article{de_rooij_conditional_2025,
	title = {Conditional universal differential equations capture population dynamics and interindividual variation in {C}-peptide production},
	volume = {11},
	doi = {10.1038/s41540-025-00570-6},
	number = {1},
	journal = {npj Systems Biology and Applications},
	author = {de Rooij, Max and van Riel, Natal A. W. and O'Donovan, Shauna D.},
	year = {2025},
	pages = {84},
}

@article{korsbo_deepnlme_2026,
	title = {{DeepNLME}: Conditional generative scientific machine learning with applications in pharmacology},
	volume = {1},
	doi = {10.70534/qm-2026-0014},
	number = {3},
	journal = {Quantitative Medicine},
	author = {Korsbo, Niklas and Tarek, Mohamed and Contento, Lorenzo and Linkevicius, Domas and Rackauckas, Christopher and Ivaturi, Vijay},
	year = {2026},
}

@article{kraft_towards_2022,
	title = {Towards hybrid modeling of the global hydrological cycle},
	volume = {26},
	doi = {10.5194/hess-26-1579-2022},
	number = {6},
	journal = {Hydrology and Earth System Sciences},
	author = {Kraft, Basil and Jung, Martin and K{\"o}rner, Marco and Koirala, Sujan and Reichstein, Markus},
	year = {2022},
	pages = {1579--1614},
}

@article{lu_neural-ode_2021,
	title = {Neural-{ODE} for pharmacokinetics modeling and its advantage to alternative machine learning models in predicting new dosing regimens},
	volume = {24},
	doi = {10.1016/j.isci.2021.102804},
	number = {7},
	journal = {iScience},
	author = {Lu, James and Deng, Kaiwen and Zhang, Xinyuan and Liu, Gengbo and Guan, Yuanfang},
	year = {2021},
	pages = {102804},
}

@article{quaghebeur_hybrid_2022,
	title = {Hybrid differential equations: Integrating mechanistic and data-driven techniques for modelling of water systems},
	volume = {213},
	doi = {10.1016/j.watres.2022.118166},
	journal = {Water Research},
	author = {Quaghebeur, Ward and Torfs, Elena and De Baets, Bernard and Nopens, Ingmar},
	year = {2022},
	pages = {118166},
}

@article{rodriguez_einns_2023,
	title = {{EINNs}: Epidemiologically-informed neural networks},
	volume = {37},
	doi = {10.1609/aaai.v37i12.26690},
	number = {12},
	journal = {Proceedings of the AAAI Conference on Artificial Intelligence},
	author = {Rodr{\'i}guez, Alexander and Cui, Jiaming and Ramakrishnan, Naren and Adhikari, Bijaya and Prakash, B. Aditya},
	year = {2023},
	pages = {14453--14460},
}

@article{schmid_wastewater-informed_2026,
	title = {Wastewater-informed neural compartmental model for long-horizon case number estimation},
	volume = {1},
	doi = {10.1038/s44482-026-00038-2},
	number = {1},
	journal = {npj Digital Public Health},
	author = {Schmid, Nina and Zacharias, Nicole and H{\"o}ser, Christoph and Bracher, Johannes and Arruda, Jonas and Papan, Cihan and Mutters, Nico T. and Hasenauer, Jan},
	year = {2026},
	pages = {35},
}

@article{philipps_current_2025,
	title = {Current state and open problems in universal differential equations for systems biology},
	volume = {11},
	doi = {10.1038/s41540-025-00550-w},
	number = {1},
	journal = {npj Systems Biology and Applications},
	author = {Philipps, Maren and Schmid, Nina and Hasenauer, Jan},
	year = {2025},
	pages = {101},
}

@article{lagergren_biologically-informed_2020,
	title = {Biologically-informed neural networks guide mechanistic modeling from sparse experimental data},
	volume = {16},
	doi = {10.1371/journal.pcbi.1008462},
	number = {12},
	journal = {PLOS Computational Biology},
	author = {Lagergren, John H. and Nardini, John T. and Baker, Ruth E. and Simpson, Matthew J. and Flores, Kevin B.},
	year = {2020},
	pages = {e1008462},
}

\appendix

\newpage

\section{Multiple indirect measurements of the predator} \label{appendix:lotka_volterra_multiple_maps}

Consider the Lotka--Volterra model from Section~\ref{section:da_lotka_volterra}, but now suppose that the prey is observed directly and the predator is measured through $k$ unknown observation maps:
\begin{align}
    \dot{N} &= rN-aNP,                        \label{eq:lotka_volterra_multiple_maps_prey} \\
    \dot{P} &= NP-mP,                         \label{eq:lotka_volterra_multiple_maps_predator} \\
    Y_0       &= N,                             \label{eq:lotka_volterra_multiple_maps_direct} \\
    Y_j     &= h_j(P), \qquad j=1,\ldots,k. \label{eq:lotka_volterra_multiple_maps_indirect}
\end{align}
Here $N>0$ and $P>0$, the parameters $r$, $a$ and $m$ are positive, and each $h_j$ has a nowhere-zero derivative. The following proposition is proved in Appendix~\ref{appendix:proof_lotka_volterra_multiple_maps}.

\begin{proposition} \label{prop:lotka_volterra_multiple_maps}
Consider the model in Equations~\eqref{eq:lotka_volterra_multiple_maps_prey}--\eqref{eq:lotka_volterra_multiple_maps_indirect}. The parameters $r$ and $m$ are structurally identifiable, whereas the parameter $a$ and every observation map $h_j$ are structurally non-identifiable. This holds no matter what the number $k$ of indirect measurements is. Specifically, given any admissible choice $(\vec{h}_1,a_1)$, every choice $(\vec{h}_2,a_2)$ of the form
\begin{align}
    a_2        &= a_1/c,                                           \notag \\
    h_{j,2}(x) &= h_{j,1}(x/c), \qquad j=1,\ldots,k,\quad c>0, \label{eq:lotka_volterra_multiple_maps_nonident}
\end{align}
gives the same observed outputs.
\end{proposition}

Thus, even observing the prey directly and measuring the predator through arbitrarily many indirect channels cannot identify the predator's absolute scale or the attack rate. This fundamentally limits what can be learnt through indirect measurements in this system: adding further uncalibrated observation maps does not remove the ambiguity.

\newpage

\section{Additional proofs}

\subsection{Proofs for Section~\ref{section:1var_models}} \label{appendix:proofs_single_variable_models}

\subsubsection{A change of state variable} \label{appendix:one_variable_change_of_state}

The five proofs of this section all rest on the same observation: two choices of the unknowns give the same observed outputs exactly when one is carried onto the other by a change of the state variable. The production--decay models of Section~\ref{section:1var_models_proddeg} are the special cases $F(X;\vec{p})=p-dX$ of Equations~\eqref{eq:1var_models_general_dynamics}--\eqref{eq:1var_models_general_obs}, so the lemma covers them too.

\begin{lemma} \label{lemma:one_variable_change_of_state}
Consider the model in Equations~\eqref{eq:1var_models_general_dynamics}--\eqref{eq:1var_models_general_obs}. Two choices $(h_1,\vec{p}_1)$ and $(h_2,\vec{p}_2)$ give the same observed outputs if and only if there is an invertible map $\sigma$ of $\mathcal{D}$ onto itself, with $\sigma'$ nowhere zero, such that
\begin{align}
    \sigma'(x)F(x;\vec{p}_1) &= F(\sigma(x);\vec{p}_2), \label{eq:one_variable_change_of_state_dynamics} \\
    h_2(\sigma(x)) &= h_1(x). \label{eq:one_variable_change_of_state_observation}
\end{align}
\end{lemma}

\begin{proof}
Suppose first that the two choices give the same observed outputs. Each observation map is differentiable with nowhere-zero derivative on the interval $\mathcal{D}$, hence strictly monotone, and so maps $\mathcal{D}$ bijectively onto its image with a differentiable inverse. Since the initial value of the output is $h_i(X(0))$ with $X(0)$ free in $\mathcal{D}$, those images are the sets of initial output values of the two choices, and therefore coincide. Hence
\begin{equation*}
    \sigma = h_2^{-1}\circ h_1
\end{equation*}
is a well-defined invertible map of $\mathcal{D}$ onto itself. Applying $h_2$ to its definition gives Equation~\eqref{eq:one_variable_change_of_state_observation}, and differentiating that relation, and dividing by $h_2'(\sigma(x))$, which is non-zero by assumption, gives
\begin{align*}
    h_2'(\sigma(x))\,\sigma'(x) &= h_1'(x) \Leftrightarrow \\
    \sigma'(x) &= \frac{h_1'(x)}{h_2'(\sigma(x))} \neq 0 .
\end{align*}

For Equation~\eqref{eq:one_variable_change_of_state_dynamics}, fix $x\in\mathcal{D}$ and let $X(t)$ be the solution of the first model with $X(0)=x$. Its output is also an output of the second model, so there is a solution $Z(t)$ of $\dot{Z}=F(Z;\vec{p}_2)$ with $h_2(Z(t))=h_1(X(t))$ for all $t$; applying $h_2^{-1}$ to both sides gives $Z(t)=\sigma(X(t))$. Differentiating this identity and using the two differential equations gives
\begin{align*}
    \dot{Z}(t) &= \sigma'(X(t))\,\dot{X}(t) = \sigma'(X(t))\,F(X(t);\vec{p}_1), \\
    \dot{Z}(t) &= F(Z(t);\vec{p}_2) = F(\sigma(X(t));\vec{p}_2).
\end{align*}
The two right-hand sides agree, and evaluating them at $t=0$, where $X(0)=x$, gives Equation~\eqref{eq:one_variable_change_of_state_dynamics}. As $x$ was arbitrary, the relation holds throughout $\mathcal{D}$.

Conversely, suppose the two relations hold, let $X(t)$ solve the first model, and put $Z(t)=\sigma(X(t))$. Then
\begin{equation*}
    \dot{Z}(t) = \sigma'(X(t))\,F(X(t);\vec{p}_1) = F(\sigma(X(t));\vec{p}_2) = F(Z(t);\vec{p}_2),
\end{equation*}
so $Z$ solves the second model, and its initial value $\sigma(X(0))$ ranges over the whole of $\mathcal{D}$ as $X(0)$ does. Equation~\eqref{eq:one_variable_change_of_state_observation} gives $h_2(Z(t))=h_1(X(t))$, so every observed output of the first choice is one of the second. For the reverse inclusion, write $y=\sigma(x)$ and note that differentiating $\sigma^{-1}(\sigma(x))=x$ gives $(\sigma^{-1})'(y)\,\sigma'(x)=1$, so that
\begin{align*}
    (\sigma^{-1})'(y)\,F(y;\vec{p}_2) &= \frac{F(\sigma(x);\vec{p}_2)}{\sigma'(x)} = F(x;\vec{p}_1) = F(\sigma^{-1}(y);\vec{p}_1), \\
    h_1(\sigma^{-1}(y)) &= h_1(x) = h_2(\sigma(x)) = h_2(y).
\end{align*}
These are the two relations with the roles of the choices exchanged and $\sigma^{-1}$ in place of $\sigma$, so the previous paragraph applies again and gives the reverse inclusion.
\end{proof}

\subsubsection{Proof of Proposition~\ref{prop:prod_deg_decay}} \label{appendix:proof_prod_deg_decay}

By Lemma~\ref{lemma:one_variable_change_of_state}, an equivalent choice satisfies
\begin{equation*}
    \sigma'(x)(-d_1x)=-d_2\sigma(x).
\end{equation*}
Since $x$ and $\sigma(x)$ are positive, integration gives $\sigma(x)=kx^{d_2/d_1}$ for some $k>0$. Hence $h_2=h_1\circ\sigma^{-1}$ gives Equation~\eqref{eq:prod_deg_decay_nonident}. Conversely, every $d_2>0$ and $k>0$ gives an admissible $\sigma$, so both $d$ and $h$ are structurally non-identifiable.

\subsubsection{Proof of Proposition~\ref{prop:prod_deg_known_production}} \label{appendix:proof_prod_deg_known_production}

By Lemma~\ref{lemma:one_variable_change_of_state}, an equivalent choice satisfies
\begin{equation*}
    \sigma'(x)(1-d_1x)=1-d_2\sigma(x).
\end{equation*}
Setting $x=1/d_1$ gives $\sigma(1/d_1)=1/d_2$. Differentiating at $x=1/d_1$ and cancelling $\sigma'(1/d_1)$ gives $d_2=d_1$. The relation then implies that $\sigma(x)=1/d_1+c(x-1/d_1)$ for some non-zero constant $c$. As $\sigma$ maps $(0,\infty)$ onto itself, it follows that $c>0$ and $\sigma(0)=0$, and hence $c=1$. Thus $\sigma(x) = x$, which implies that $h_2=h_1$, so both $d$ and $h$ are globally structurally identifiable.

\subsubsection{Proof of Proposition~\ref{prop:prod_deg_unknown_production}} \label{appendix:proof_prod_deg_unknown_production}

By Lemma~\ref{lemma:one_variable_change_of_state},
\begin{equation*}
    \sigma'(x)(p_1-d_1x)=p_2-d_2\sigma(x).
\end{equation*}
Differentiating the original relation with respect to $x$ gives
\begin{equation*}
    \sigma''(x)(p_1-d_1x)-d_1\sigma'(x)=-d_2\sigma'(x).
\end{equation*}
Evaluating this at $x=p_1/d_1$ gives
\begin{equation*}
    -d_1\sigma'\!\left(\frac{p_1}{d_1}\right)
    =-d_2\sigma'\!\left(\frac{p_1}{d_1}\right).
\end{equation*}
Because $\sigma'$ is nowhere zero, we can cancel it and obtain $d_1=d_2$. Write their common value as $d$. The original relation can then be rearranged as
\begin{equation*}
    \sigma'(x)\left(x-\frac{p_1}{d}\right)
    =\sigma(x)-\frac{p_2}{d}.
\end{equation*}
Integrating gives
\begin{equation*}
    \sigma(x)=\frac{p_2}{d}+c\left(x-\frac{p_1}{d}\right)
\end{equation*}
for some non-zero constant $c$. Since $\sigma$ maps $(0,\infty)$ onto itself, we must have $c>0$ and $\sigma(0)=0$. Therefore
\begin{equation*}
    0=\frac{p_2}{d}-c\frac{p_1}{d},
    \qquad\text{so}\qquad
    c=\frac{p_2}{p_1}.
\end{equation*}
It follows that
\begin{equation*}
    \sigma(x)=\frac{p_2}{p_1}x,
    \qquad
    \sigma^{-1}(x)=\frac{p_1}{p_2}x.
\end{equation*}
Hence
\begin{equation*}
    h_2(x)=h_1\!\left(\sigma^{-1}(x)\right)
          =h_1\!\left(\frac{p_1x}{p_2}\right),
\end{equation*}
which is Equation~\eqref{eq:prod_deg_unknown_production_nonident}. Conversely, every $p_2>0$ in this expression gives an equivalent choice.

\subsubsection{Proof of Proposition~\ref{prop:1var_general_nonidentifiability}} \label{appendix:proof_1var_general_nonidentifiability}

Fix two admissible parameter vectors $\vec{p}_1$ and $\vec{p}_2$, and pick an arbitrary reference state $x_i^0\in\mathcal{D}$ for each. Define
\begin{equation*}
    \tau_i(x)=\int_{x_i^0}^{x}\frac{\mathrm{d}s}{F(s;\vec{p}_i)},
\end{equation*}
so that $\tau_i(x)$ is the time at which the $i$th model is in state $x$, its clock being set to zero in state $x_i^0$. Because $F$ has no steady state and every solution exists for all time, each $\tau_i$ maps $\mathcal{D}$ bijectively onto $\mathbb{R}$, and so is a time coordinate on the state space. Therefore $\sigma=\tau_2^{-1}\circ\tau_1$, which matches the states the two models occupy at the same time, is an admissible change of state variable and satisfies Equation~\eqref{eq:one_variable_change_of_state_dynamics}. Lemma~\ref{lemma:one_variable_change_of_state} then gives an equivalent choice for every $\vec{p}_2$. Taking $\vec{p}_2=\vec{p}_1$, the maps $\sigma_c=\tau_1^{-1}\circ(\tau_1+c)$ for non-zero $c$ are non-identity admissible changes of state variable and hence give distinct equivalent observation maps. Thus every parameter component that can vary, and $h$, is structurally non-identifiable.

\subsubsection{Proof of Proposition~\ref{prop:1var_general_identifiability}} \label{appendix:proof_1var_general_identifiability}

The hypothesis that the steady states and their eigenvalues determine $\vec{p}$ presupposes that there are finitely many of them and that each is simple.

Fix an admissible choice $(h_1,\vec{p}_1)$ and let $(h_2,\vec{p}_2)$ give the same observed outputs. Lemma~\ref{lemma:one_variable_change_of_state} supplies a change of state variable $\sigma$ with
\begin{equation}
    \sigma'(x)F(x;\vec{p}_1)=F(\sigma(x);\vec{p}_2), \label{eq:1var_general_proof_relation}
\end{equation}
and $\sigma$ is increasing because both observation maps have the same direction.

An increasing bijection of $\mathcal{D}$ onto itself carries the zeros of $F(\cdot;\vec{p}_1)$ onto the zeros of $F(\cdot;\vec{p}_2)$ and pairs them off in increasing order. Differentiating Equation~\eqref{eq:1var_general_proof_relation} at a zero $x^*$ of $F(\cdot;\vec{p}_1)$, where the term carrying $\sigma''$ has the factor $F(x^*;\vec{p}_1)=0$ and drops out, leaves
\begin{equation*}
    \sigma'(x^*)F'(x^*;\vec{p}_1)=F'(\sigma(x^*);\vec{p}_2)\,\sigma'(x^*).
\end{equation*}
Paired steady states therefore have equal eigenvalues, and the assumption of the proposition gives $\vec{p}_2=\vec{p}_1$. We drop the parameters from the notation and write Equation~\eqref{eq:1var_general_proof_relation} as $\sigma'(x)F(x)=F(\sigma(x))$.

The two models are now one and the same, so $\sigma$ maps the steady states onto themselves in increasing order and hence fixes each of them; it therefore maps each interval between two adjacent zeros of $F$ onto itself. Let $J$ be such an interval and $\tau(x)=\int_{x_J}^{x}\mathrm{d}s/F(s)$ for a base point $x_J\in J$, which is defined because $F$ has no zero on $J$. Then
\begin{equation*}
    \frac{\mathrm{d}}{\mathrm{d}x}\tau(\sigma(x))=\frac{\sigma'(x)}{F(\sigma(x))}=\frac{1}{F(x)}=\tau'(x),
    \qquad\text{so}\qquad
    \tau(\sigma(x))=\tau(x)+c_J
\end{equation*}
for a constant $c_J$: on $J$, the map $\sigma$ moves every state forward along the flow by the time $c_J$.

Let $J$ now be an interval carrying a solution that reaches an end of $\mathcal{D}$ in finite time, say $X_{\max}$; the remaining cases are the same with the two ends of $\mathcal{D}$, or the direction of time, exchanged. The time that the solution through $x\in J$ still needs in order to leave $\mathcal{D}$ is
\begin{equation*}
    R(x)=\int_{x}^{X_{\max}}\frac{\mathrm{d}s}{F(s)}=\lim_{z\to X_{\max}}\tau(z)-\tau(x),
\end{equation*}
which by hypothesis is finite, is positive, and is arbitrarily small for $x$ close to $X_{\max}$. The previous display gives $R(\sigma(x))=R(x)-c_J$, so a positive $c_J$ would make $R$ negative somewhere on $J$; applying the same to $\sigma^{-1}$, which moves every state forward by $-c_J$, excludes a negative $c_J$. Hence $c_J=0$ and $\sigma$ is the identity on $J$.

Let $x^*$ be an end of $J$ that is a steady state, and $J'$ the interval on its other side. Moving forward by a time $c$ fixes $x^*$ and has the derivative $\exp(F'(x^*)c)$ there, so the one-sided derivatives of $\sigma$ at $x^*$ are $1$ from $J$ and $\exp(F'(x^*)c_{J'})$ from $J'$. These agree because $\sigma$ is differentiable at $x^*$, and $F'(x^*)\neq0$ then gives $c_{J'}=0$. As there are finitely many steady states, repeating this across each of them gives $\sigma(x)=x$ throughout $\mathcal{D}$. Equation~\eqref{eq:one_variable_change_of_state_observation} then gives $h_2=h_1$, so both $\vec{p}$ and $h$ are globally structurally identifiable.

\subsection{Proofs for Section~\ref{section:diff_algebra}} \label{appendix:proofs_differential_algebra}

\subsubsection{Proof of identifiability of the decay chain with a calibration sample} \label{appendix:proof_chain_calibration}

Suppose that a calibration sample gives the value $y_0=h(x_0)$ of the observation map at some state $x_0>0$. By Proposition~\ref{prop:chain}, any choice giving the same observed outputs as $(h_1,\lambda_1,\lambda_2)$ has the same rates and an observation map $h_2(x)=h_1(x/k)$ with $k>0$. Both observation maps must also match the calibration sample, so
\begin{equation*}
    h_1(x_0/k)=h_2(x_0)=y_0=h_1(x_0).
\end{equation*}
Since $h_1$ has a nowhere-zero derivative, it is injective, which gives $x_0/k=x_0$ and hence $k=1$. Thus $h_2=h_1$, and $h$, $\lambda_1$ and $\lambda_2$ are all globally structurally identifiable.

\subsubsection{Proof of Proposition~\ref{prop:sir}} \label{appendix:proof_sir}

The variable $R$ is unobserved and does not appear in Equations~\eqref{eq:sir_dynamics_s} and~\eqref{eq:sir_dynamics_i}, so we eliminate $S$ using these two equations alone. First, Equation~\eqref{eq:sir_dynamics_i} gives
\begin{equation}
    S=\frac{\dot{I}+\gamma I}{\beta I}. \label{eq:sir_proof_s}
\end{equation}
Differentiating this expression gives
\begin{equation}
    \dot{S}=\frac{I\ddot{I}-\dot{I}^2}{\beta I^2}. \label{eq:sir_proof_s_derivative}
\end{equation}
Substituting Equations~\eqref{eq:sir_proof_s} and~\eqref{eq:sir_proof_s_derivative} into $\dot{S}=-\beta SI$ yields
\begin{align}
    \frac{I\ddot{I}-\dot{I}^2}{\beta I^2}
    &=-\beta\left(\frac{\dot{I}+\gamma I}{\beta I}\right)I \Leftrightarrow \notag \\
     \ddot{I} &=\frac{\dot{I}^2}{I}-\beta I\dot{I}-\beta\gamma I^2. \label{eq:sir_proof_eliminated}
\end{align}
Finally, write $g=h^{-1}$, which gives $I=g(Y)$. Differentiating once and twice gives
\begin{align*}
    \dot{I}=g'(Y)\dot{Y}, \\
    \ddot{I}=g''(Y)\dot{Y}^2+g'(Y)\ddot{Y}.
\end{align*}
Substituting these expressions into Equation~\eqref{eq:sir_proof_eliminated} and dividing by $g'(Y)$ gives the input--output equation
\begin{equation}
    \ddot{Y}
    =\left(\frac{g'(Y)}{g(Y)}-\frac{g''(Y)}{g'(Y)}\right)\dot{Y}^2
     -\beta g(Y)\dot{Y}
     -\beta\gamma\frac{g(Y)^2}{g'(Y)}. \label{eq:sir_proof_io}
\end{equation}

\medskip\noindent
The coefficients of the input--output equation are the identifiable quantities. Suppose that two choices give the same observed outputs, and write $g_i=h_i^{-1}$. Equating the coefficients of $\dot{Y}$ gives
\begin{equation}
    \beta_1g_1(Y)=\beta_2g_2(Y). \label{eq:sir_proof_g_relation}
\end{equation}
Thus $g_2=cg_1$, where $c=\beta_1/\beta_2>0$. This relation also makes the coefficients of $\dot{Y}^2$ equal. Equality of the remaining coefficients then gives
\begin{equation*}
    \beta_1\gamma_1\frac{g_1^2}{g_1'}
    =\beta_2\gamma_2\frac{g_2^2}{g_2'}
    =\beta_1\gamma_2\frac{g_1^2}{g_1'},
\end{equation*}
so $\gamma_2=\gamma_1$. Since $g_2=cg_1$ is equivalent to $h_2(x)=h_1(x/c)$, any two indistinguishable choices must have
\begin{equation*}
    \beta_2=\frac{\beta_1}{c}, \qquad
    \gamma_2=\gamma_1, \qquad
    h_2(x)=h_1(x/c), \qquad c>0.
\end{equation*}
Conversely, for any $c>0$, scaling all three states by $c$ takes a solution for $(\beta_1,\gamma_1)$ to a solution for $(\beta_1/c,\gamma_1)$, while $h_2(cI)=h_1(I)$. These choices therefore give exactly the same observed outputs, proving Equation~\eqref{eq:sir_nonident}. Hence $\gamma$ is globally structurally identifiable, whereas $\beta$ and $h$ are not.

If instead the total population size $N=S+I+R$ is known, then $I$ ranges over $(0,N)$. Each inverse observation map $g_i$ must therefore have range $(0,N)$. Since $g_2=cg_1$ has range $(0,cN)$, this is possible only when $c=1$. It follows that $h_2=h_1$, $\beta_2=\beta_1$ and $\gamma_2=\gamma_1$, so all three unknowns are globally structurally identifiable.

\subsubsection{Proof of Proposition~\ref{prop:lotka_volterra}} \label{appendix:proof_lotka_volterra}

First, Equation~\eqref{eq:lotka_volterra_dynamics_pred} gives
\begin{equation}
    N=\frac{\dot{P}+mP}{P}. \label{eq:lotka_volterra_proof_n}
\end{equation}
Differentiating this expression gives
\begin{equation}
    \dot{N}=\frac{P\ddot{P}-\dot{P}^2}{P^2}. \label{eq:lotka_volterra_proof_n_derivative}
\end{equation}
Substituting Equations~\eqref{eq:lotka_volterra_proof_n} and~\eqref{eq:lotka_volterra_proof_n_derivative} into $\dot{N}=rN-aNP$ yields
\begin{align}
    \frac{P\ddot{P}-\dot{P}^2}{P^2}
    &=\left(\frac{\dot{P}+mP}{P}\right)(r-aP) \Leftrightarrow \notag \\
    \ddot{P} &=\frac{\dot{P}^2}{P}+(r-aP)\dot{P}+mP(r-aP). \label{eq:lotka_volterra_proof_eliminated}
\end{align}
Finally, write $g=h^{-1}$, which gives $P=g(Y)$. Differentiating once and twice gives
\begin{align*}
    \dot{P}=g'(Y)\dot{Y}, \\
    \ddot{P}=g''(Y)\dot{Y}^2+g'(Y)\ddot{Y}.
\end{align*}
Substituting these expressions into Equation~\eqref{eq:lotka_volterra_proof_eliminated} and dividing by $g'(Y)$ gives the input--output equation
\begin{equation}
    \ddot{Y}
    =\left(\frac{g'(Y)}{g(Y)}-\frac{g''(Y)}{g'(Y)}\right)\dot{Y}^2
     +\bigl(r-ag(Y)\bigr)\dot{Y}
     +m\,\frac{g(Y)\bigl(r-ag(Y)\bigr)}{g'(Y)}. \label{eq:lotka_volterra_proof_io}
\end{equation}

\medskip\noindent
The coefficients of the input--output equation are the identifiable quantities. Suppose that two choices give the same observed outputs, and write $g_i=h_i^{-1}$. Equating the coefficients of $\dot{Y}$ gives
\begin{equation}
    r_1-a_1g_1(Y)=r_2-a_2g_2(Y), \label{eq:lotka_volterra_proof_g_relation}
\end{equation}
that is, $g_2=cg_1+b$, where $c=a_1/a_2>0$ and $b=(r_2-r_1)/a_2$. Since $P$ takes every value in $(0,\infty)$, both $g_1$ and $g_2$ have range $(0,\infty)$, whereas the range of $cg_1+b$ is $(b,\infty)$. Hence $b=0$, so
\begin{equation*}
    g_2=cg_1, \qquad r_2=r_1, \qquad a_2=a_1/c, \qquad c>0.
\end{equation*}
This relation also makes the coefficients of $\dot{Y}^2$ equal. It also gives $g_2'=cg_1'$ and $a_2g_2=a_1g_1$, so equality of the remaining coefficients reads
\begin{equation*}
    m_1\frac{g_1(r_1-a_1g_1)}{g_1'}
    =m_2\frac{g_2(r_2-a_2g_2)}{g_2'}
    =m_2\frac{g_1(r_1-a_1g_1)}{g_1'},
\end{equation*}
where the shared factor vanishes only where $g_1(Y)=r_1/a_1$, which is at most a single point because $g_1$ is injective. Hence $m_2=m_1$. Since $g_2=cg_1$ is equivalent to $h_2(x)=h_1(x/c)$, any two indistinguishable choices must have
\begin{equation*}
    r_2=r_1, \qquad
    a_2=\frac{a_1}{c}, \qquad
    m_2=m_1, \qquad
    h_2(x)=h_1(x/c), \qquad c>0.
\end{equation*}
Conversely, for any $c>0$, scaling the predator by $c$ takes a solution for $(r_1,a_1,m_1)$ to a solution for $(r_1,a_1/c,m_1)$, while $h_2(cP)=h_1(P)$. These choices therefore give exactly the same observed outputs, proving Equation~\eqref{eq:lotka_volterra_nonident}. Hence $r$ and $m$ are globally structurally identifiable, whereas $a$ and $h$ are not.

\subsection{Proofs for Section~\ref{section:reporters}} \label{appendix:proofs_reporters}

\subsubsection{Proof of Proposition~\ref{prop:reporter_reduction}} \label{appendix:proof_reporter_reduction}

Integrating Equation~\eqref{eq:reporter_marker}, with $d_j$ known, gives
\begin{equation*}
    Y_j(t) = R_j(t) = e^{-d_jt}R_j(0)+\int_0^t e^{-d_j(t-s)}f_j\bigl(\vec{u}_j(s)\bigr)\,\mathrm{d}s.
\end{equation*}
The reporter output is therefore determined by the output $f_j\bigl(\vec{u}_j(t)\bigr)$ of the model in Equation~\eqref{eq:reporter_reduced}, along the same solution of Equation~\eqref{eq:reporter_ode}, together with the initial reporter level $R_j(0)$. Conversely, differentiating this expression gives
\begin{equation*}
    \dot{Y}_j(t)+d_jY_j(t) = f_j\bigl(\vec{u}_j(t)\bigr),
\end{equation*}
so the reporter output determines the output of the reduced model. For every initial value of $\vec{X}$, the same initial reporter levels are admissible, whatever the unknowns. Therefore, two choices of $(\vec{f},\vec{p})$ give the same reporter outputs if and only if they give the same outputs for the reduced model.

\subsubsection{Proof of Proposition~\ref{prop:reporter_decay}} \label{appendix:proof_reporter_decay}

\begin{proof}
Let $(\vec{f}_1,\vec{p}_1,\vec{d}_1)$ and $(\vec{f}_2,\vec{p}_2,\vec{d}_2)$ give the same observed outputs, and fix $j$. Fix an admissible initial value $\vec{x}_1$, and let $\vec{X}_1(t)$ be the resulting state trajectory for the first choice, so $\vec{X}_1(0)=\vec{x}_1$. Let $\vec{u}_{j,1}$ be the relevant sub-vector of $\vec{x}_1$, and write $q=f_{j,1}(\vec{u}_{j,1})\geq0$. Keeping $\vec{x}_1$ fixed, choose any $R_{j,1}(0)=r\geq0$, and write $\vec{Y}^{(r)}$ for the resulting output. Then
\begin{equation*}
    Y_j^{(r)}(0)=r, \qquad \dot{Y}_j^{(r)}(0)=q-d_{j,1}r.
\end{equation*}
The second choice has a solution with this same output $\vec{Y}^{(r)}$, from some initial value $\vec{x}_2(r)$. Let $\vec{u}_{j,2}(r)$ be the relevant sub-vector of $\vec{x}_2(r)$. Its reporter equation gives
\begin{equation*}
    \dot{Y}_j^{(r)}(0)=f_{j,2}(\vec{u}_{j,2}(r))-d_{j,2}Y_j^{(r)}(0).
\end{equation*}
Combining the two displayed equations gives
\begin{equation*}
    f_{j,2}(\vec{u}_{j,2}(r))=q+(d_{j,2}-d_{j,1})r.
\end{equation*}
The initial value $\vec{x}_2(r)$ may vary with $r$, but the left-hand side is always non-negative, whereas $q$ is fixed. Since $r$ can be arbitrarily large, this is impossible if $d_{j,2}<d_{j,1}$. Thus $d_{j,2}\geq d_{j,1}$. The same argument with the two choices exchanged gives $d_{j,1}\geq d_{j,2}$, so $d_{j,1}=d_{j,2}$. Since $j$ was arbitrary, every $d_j$ is globally structurally identifiable.
\end{proof}

\subsubsection{Proof of Proposition~\ref{prop:reporter_decay_signal}} \label{appendix:proof_reporter_decay_signal}

Since $f_1$ is positive, Proposition~\ref{prop:reporter_decay} gives $d_2=d_1$. With the decay rate known, Proposition~\ref{prop:reporter_reduction} says that two choices of the remaining unknowns give the same observed outputs if and only if they do so for the observation-map model 
\begin{align*}
    \dot{X} &= -aX, \\
    Y &=f(X).
\end{align*}
This is the model of Proposition~\ref{prop:prod_deg_decay}, with $a$ in place of $d$ and $f$ in place of $h$, and Equation~\eqref{eq:prod_deg_decay_nonident} gives Equation~\eqref{eq:reporter_decay_nonident}.

\subsubsection{Proof of Proposition~\ref{prop:reporter_goodwin}} \label{appendix:proof_reporter_goodwin}

Since $f$ is positive, Proposition~\ref{prop:reporter_decay} identifies $d$, and Proposition~\ref{prop:reporter_reduction} reduces the remaining problem to the observation-map model formed by Equations~\eqref{eq:reporter_goodwin_x1}--\eqref{eq:reporter_goodwin_x3} with $Y=f(X_3)$.

First, write $z=X_3$. Equations~\eqref{eq:reporter_goodwin_x2} and~\eqref{eq:reporter_goodwin_x3} give
\begin{align}
    X_2 &= \dot{z}+b_3z, \label{eq:goodwin_proof_x2} \\
    X_1 &= \ddot{z}+(b_2+b_3)\dot{z}+b_2b_3z. \label{eq:goodwin_proof_x1}
\end{align}
Differentiating Equation~\eqref{eq:goodwin_proof_x1} and substituting it into Equation~\eqref{eq:reporter_goodwin_x1} gives
\begin{equation}
    \dddot{z}=-\alpha_1\ddot{z}-\alpha_2\dot{z}-\alpha_3z+\frac{a}{1+z^m}, \label{eq:goodwin_proof_eliminated}
\end{equation}
where
\begin{equation*}
    \alpha_1=b_1+b_2+b_3, \qquad
    \alpha_2=b_1b_2+b_1b_3+b_2b_3, \qquad
    \alpha_3=b_1b_2b_3.
\end{equation*}

Finally, write $g=f^{-1}$, which gives $z=g(Y)$. Differentiation gives
\begin{align*}
    \dot{z}   &=g'(Y)\dot{Y}, \\
    \ddot{z}  &=g''(Y)\dot{Y}^2+g'(Y)\ddot{Y}, \\
    \dddot{z} &=g'''(Y)\dot{Y}^3+3g''(Y)\dot{Y}\ddot{Y}+g'(Y)\dddot{Y}.
\end{align*}
Substituting these expressions into Equation~\eqref{eq:goodwin_proof_eliminated} and dividing by $g'(Y)$ gives the input--output equation
\begin{align}
    \dddot{Y}
    ={}&-3\frac{g''(Y)}{g'(Y)}\dot{Y}\ddot{Y}
        -\frac{g'''(Y)}{g'(Y)}\dot{Y}^3
        -\alpha_1\ddot{Y}
        -\alpha_1\frac{g''(Y)}{g'(Y)}\dot{Y}^2 \notag \\
       &-\alpha_2\dot{Y}
        -\alpha_3\frac{g(Y)}{g'(Y)}
        +\frac{a}{g'(Y)\bigl(1+g(Y)^m\bigr)}. \label{eq:goodwin_proof_io}
\end{align}

The coefficients of the input--output equation are the identifiable quantities. Suppose that two choices $g_1$ and $g_2$ give the same observed dynamics. Equating the coefficients of $\dot{Y}\ddot{Y}$ gives
\begin{equation*}
    \frac{g_1''(Y)}{g_1'(Y)}=\frac{g_2''(Y)}{g_2'(Y)},
\end{equation*}
and hence $g_2=cg_1+e$ for constants $c\neq0$ and $e$. Both inverse maps have range $(0,\infty)$, so $c>0$ and $e=0$. Equating the coefficients of $\ddot{Y}$ and $\dot{Y}$ also shows that $\alpha_1$ and $\alpha_2$ are identifiable.

It remains to compare the constant coefficients. Setting $x=g_1(Y)$ and using $g_2=cg_1$ gives
\begin{equation}
    (\alpha_{3,2}-\alpha_{3,1})x+\frac{a_1}{1+x^{m_1}}
    =\frac{a_2/c}{1+(cx)^{m_2}}. \label{eq:goodwin_proof_constant}
\end{equation}
Letting $x\to\infty$ shows that $\alpha_{3,2}=\alpha_{3,1}$, and letting $x\to0$ then gives $a_2=ca_1$. Equation~\eqref{eq:goodwin_proof_constant} therefore leaves
\begin{equation*}
    x^{m_1}=c^{m_2}x^{m_2}, \qquad x>0,
\end{equation*}
so $m_2=m_1$ and $c=1$. Thus $a_2=a_1$, $g_2=g_1$, and hence $f_2=f_1$.

The input--output equation identifies the three decay rates as an unordered set. To determine their order, use the reconstruction equations above as constraints. For each fixed $z>0$, Equation~\eqref{eq:goodwin_proof_x2} and the condition $X_2>0$ show that the possible values of $\dot{z}$ are exactly
\begin{equation*}
    \dot{z}>-b_3z.
\end{equation*}
The boundary of this observed range identifies $b_3$. At fixed $z$ and $\dot{z}=0$, Equation~\eqref{eq:goodwin_proof_x1} similarly shows that the possible values of $\ddot{z}$ are exactly
\begin{equation*}
    \ddot{z}>-b_2b_3z.
\end{equation*}
Its boundary identifies $b_2$, since $b_3$ is already known, and $\alpha_1=b_1+b_2+b_3$ then identifies $b_1$. Therefore all the unknowns are globally structurally identifiable.

\subsection{Proofs for Appendix~\ref{appendix:lotka_volterra_multiple_maps}} \label{appendix:proofs_lotka_volterra_multiple_maps}

\subsubsection{Proof of Proposition~\ref{prop:lotka_volterra_multiple_maps}} \label{appendix:proof_lotka_volterra_multiple_maps}

Let $(N(t),P(t))$ solve the model with $(\vec{h}_1,r,a_1,m)$, fix $c>0$, and put $\widetilde{N}=N$ and $\widetilde{P}=cP$. Then
\begin{align*}
    \dot{\widetilde{N}} &= r\widetilde{N}-(a_1/c)\widetilde{N}\widetilde{P}, \\
    \dot{\widetilde{P}} &= \widetilde{N}\widetilde{P}-m\widetilde{P},
\end{align*}
so $(\widetilde{N},\widetilde{P})$ solves the model with attack rate $a_2=a_1/c$. Its directly observed output is still $\widetilde{N}=N$, while Equation~\eqref{eq:lotka_volterra_multiple_maps_nonident} gives
\begin{equation*}
    h_{j,2}(\widetilde{P})=h_{j,1}(\widetilde{P}/c)=h_{j,1}(P), \qquad j=1,\ldots,k.
\end{equation*}
The rescaling $(N,P)\mapsto(N,cP)$ is a bijection of the positive state space, so the two choices give the same set of observed outputs. For $c\neq1$, we have $a_2\neq a_1$ and $h_{j,2}\neq h_{j,1}$ for every $j$, since each $h_{j,1}$ has a nowhere-zero derivative. Thus $a$ and all of the observation maps are structurally non-identifiable. The parameters $r$ and $m$ are structurally identifiable by Proposition~\ref{prop:lotka_volterra}, since two choices that give the same outputs here also give the same output $Y_1=h_1(P)$.

\newpage

\section{Numerical implementation} \label{appendix:numerical_implementation}

\subsection{Synthetic data generation} \label{appendix:sir_data_generation}
The data used in Figure~\ref{fig:sir_fits} are generated from Equations~\eqref{eq:sir_dynamics_s}--\eqref{eq:sir_obs} with parameter values $\beta=8\times10^{-4}$ and $\gamma=0.2$, initial condition $(S,I,R)(0)=(990,10,0)$, and observation map
\begin{equation*}
    h(x)=\frac{950x}{x+500}.
\end{equation*}
The total population size is thus $N=1000$. The model is solved on $t\in[0,40]$ using the classical fourth-order Runge--Kutta method with a fixed time step of $0.02$. We retain only the measurement $Y=h(I)$, sampled at $51$ evenly spaced time points ($\Delta t=0.8$), and treat $S$, $I$ and $R$ as latent throughout. No observation noise is added.

\subsection{Fitting procedure} \label{appendix:sir_fitting}
We represent the unknown observation map by a neural network with a single hidden layer of six $\tanh$ units,
\begin{equation*}
    h(x)=\theta_0+\sum_{k=1}^{6}\theta_{1,k}\tanh\left(\theta_{2,k}x+\theta_{3,k}\right),
\end{equation*}
where the weights $\theta_{1,k}$ and $\theta_{2,k}$ are kept positive by writing them as softplus functions, $\log(1+e^{z})$, of unconstrained parameters. This makes $h$ strictly increasing, so that it has a nowhere-zero derivative, as assumed in Proposition~\ref{prop:sir}. For numerical conditioning, the states and measurements are rescaled by a factor $10^{-3}$ before fitting, and the network acts on the rescaled $I$. This is a fixed change of units, and gives the fit no information about $N$.

We consider two cases. When $N$ is known, the unknowns are $\beta$, $\gamma$, $I(0)$, and the network parameters, with $S(0)=N-I(0)$, that is, with $R(0)=0$. When $N$ is unknown, $S(0)$ is also fitted. The positive scalar unknowns are fitted on a logarithmic scale. For each choice of the unknowns, we solve Equations~\eqref{eq:sir_dynamics_s}--\eqref{eq:sir_dynamics_i} as in Appendix~\ref{appendix:sir_data_generation}, and compute the loss
\begin{equation*}
    \mathcal{L}=\sum_{i=1}^{51}\bigl(Y(t_i)-h(I(t_i))\bigr)^2,
\end{equation*}
where $Y(t_i)$ are the synthetic measurements.

For each case, we carry out 30 independent fits from random starting points. The starting values of $\beta$, $\gamma$, and the initial conditions are drawn log-uniformly from $\beta\in[1.6\times10^{-4},4\times10^{-3}]$, $\gamma\in[0.04,1]$, and $I(0)\in[0.1,900]$ when $N$ is known, and from the same ranges for $\beta$ and $\gamma$, together with $S(0)\in[50,5000]$ and $I(0)\in[0.1,2000]$, when $N$ is unknown. The unconstrained parameters behind $\theta_{1,k}$ and $\theta_{2,k}$ are drawn from $\mathcal{N}(-2,1)$ and $\mathcal{N}(1,1)$, $\theta_{3,k}$ from $\mathcal{N}(-1,4)$, and $\theta_0$ starts at $0$. Each fit runs $20{,}000$ steps of the Adam optimiser~\cite{kingma_adam_2015} with learning rate $5\times10^{-3}$, followed by up to $5{,}000$ iterations of L-BFGS~\cite{liu_limited_1989} with a backtracking line search, as implemented in Optim.jl~\cite{mogensen_optim_2018}. Gradients are computed by forward-mode automatic differentiation using ForwardDiff.jl~\cite{revels_forward-mode_2016}. All computations are carried out in Julia~\cite{bezanson_julia_2017}.

We assess each fit by its root-mean-square error, divided by the range of the measurements, and retain fits for which this is below $3\times10^{-4}$. All 30 fits with $N$ known, and 27 of the 30 fits with $N$ unknown, meet this criterion. The remaining three stopped at a local minimum with a normalised error of $0.34$. Panel B of Figure~\ref{fig:sir_fits} shows the 20 retained fits with the lowest loss for each case, and panels C and D the 10 lowest of these. In panel B, each fitted network is evaluated across the full axis, including values of $x$ outside the range of $I$ visited by its fit, where the data do not constrain it.

Across the 20 fits shown with $N$ known, the recovered parameters lie within $\beta\in[7.994,8.000]\times10^{-4}$ and $\gamma\in[0.1998,0.2000]$, matching the true values. With $N$ unknown, $\gamma\in[0.1998,0.2001]$ is still recovered, whereas $\beta$ ranges over $[1.97,15.6]\times10^{-4}$. This corresponds to scaling factors $c$ between about $0.5$ and $4$ in Proposition~\ref{prop:sir}, which is the numerical signature of the non-identifiability of $\beta$ and $h$.

\end{document}